\documentclass[11pt]{amsart}

\usepackage{amsthm}
\usepackage{amssymb}
\usepackage{amsmath} 

\usepackage{hyperref}

\newcommand{\bea}{\begin{eqnarray*}}
\newcommand{\eea}{\end{eqnarray*}}
\newcommand{\bean}{\begin{eqnarray}}
\newcommand{\eean}{\end{eqnarray}}
\newtheorem{thm}{Theorem}[section]
\newtheorem{prop}[thm]{Proposition}
\newtheorem{lem}[thm]{Lemma}
\newtheorem{cor}[thm]{Corollary}
\newtheorem{rem}[thm]{Remark}
\newtheorem{defn}[thm]{Definition}

\newtheorem{ex}[thm]{Example}

\newcommand{\sgn}{\mathrm{sgn}}
\newcommand{\HH}{\mathcal{H}_n}
\newcommand{\wt}{\Lambda}

\begin{document}
\title[Polynomial repr. of $(C^\vee_n, C_n)$-DAHA
 for specialized parameters]
{The polynomial representation
 of the double affine Hecke algebra
 of type $(C^\vee_n, C_n)$ for specialized parameters}
\date{}
\author{Masahiro Kasatani}
\address{Department of Mathematics, Graduate School of Science,
Kyoto University, Kyoto 606-8502, Japan}
\email{kasatani@math.kyoto-u.ac.jp}

\begin{abstract}
In this paper,
we study the polynomial representation of the double affine Hecke algebra
 of type $(C^\vee_n, C_n)$ for specialized parameters.
Inductively and combinatorially,
we give a linear basis of the representation
in terms of linear combinations of non-symmetric Koornwinder polynomials.
The basis consists of generalized eigenfunctions with respect to
 $q$-Dunkl-Cherednik operators $\widehat{Y}_i$,
and it gives a way to cancel out poles of non-symmetric Koornwinder polynomials.
We examine irreducibility and $Y$-semisimplicity of the representation
for the specialized parameters.
For some cases,
 we give a characterization of the subrepresentations by vanishing conditions
for Laurent polynomials.
\end{abstract}

\maketitle

\section{Introduction}

In 1990's, Cherednik introduced the double affine Hecke algebra (DAHA).
(\cite{ChBook}).
It is a unital associative algebra with some parameters
attached to affine root systems,
and it contains two affine Hecke algebras.
The DAHA has an action on a ring of multivariable Laurent polynomials,
 which is called a polynomial representation.
For reduced affine root systems,
he proved irreducibility of the polynomial representation for generic parameters,
and solved conjectures (duality relations, evaluation formulas, and so on)
 on Macdonald polynomials\cite{Ma}
 by exploiting an anti-involution of the DAHA.

For the type $(C^\vee_n, C_n)$, which is a non-reduced affine root system,
Noumi introduced a polynomial representation
of the affine Hecke algebra of type $C_n$
 to study Macdonald-Koornwinder polynomials (\cite{No}).
His realization of operators is called the Noumi representation.
Standing on Noumi's work,
Sahi introduced the double affine Hecke algebra
 of type $(C^\vee_n, C_n)$ (\cite{Sa})
and an extension of the polynomial representation to DAHA-module.
He proved irreducibility of the representation
 for generic parameters, and the duality relations.
In \cite{St}, the evaluation formulas are shown by Stokman.

For rank $n=1$,
Oblomkov and Stoica classified finite dimensional representations
of the DAHA of type $(C^\vee_1, C_1)$ (\cite{ObSt})
in the case when a parameter $q$ is not a root of unity.
They gave a description of subrepresentations in terms of
 the polynomial representation.
The case when $q=1$ and the other parameters are generic is
studied in \cite{Ob}.

Recently, Cherednik gave a condition when the polynomial representation
of the DAHA is irreducible (\cite{Ch})
for reduced root systems and a generic parameter $q$.
He introduced a pairing on the representation,
and determined irreducibility by the condition 
if the radical of the pairing is zero or not.
He also determined whether the radical is zero or not
by so-called affine exponents.

We will recall the proof of irreducibility of the polynomial representation
for generic parameters.
$q$-Dunkl-Cherednik operators $\widehat{Y}_i$ are simultaneously diagonalizable
on the representation.
This property is called $Y$-semisimplicity.
$Y$-eigenfunctions are called non-symmetric Macdonald polynomials.
(For type $(C^\vee_n, C_n)$,
 they are also called non-symmetric Koornwinder polynomials.)
There exist operators which send a $Y$-eigenfunction to
 other $Y$-eigenfunctions.
These operators are called intertwiners.
By applying the intertwiners, we see that
 any $Y$-eigenfunction is cyclic.

In this paper,
we treat the DAHA of type $(C^\vee_n, C_n)$ with $n\geq2$.
The algebra has 6 parameters $q$, $t$, $a^*$, $b^*$, $c^*$, and $d^*$.
We study the polynomial representation for specialized parameters
 where the representation can be non-$Y$-semisimple or reducible.
(The parameter $q$ may be a root of unity.)
For such specialized parameters,
 non-symmetric Koornwinder polynomials may have some poles,
 or the argument in terms of the intertwiners above is not enough.
In order to overcome these difficulties,
we introduce certain linear combinations of non-symmetric
Koornwinder polynomials where the poles are cancelled out
(we call them {\it modified polynomials}).
Modified polynomials are generalized $Y$-eigenfunctions.
We also introduce {\it modified intertwiners}.
Using them,
we give an inductive and combinatorial method
to compute modified polynomials.
So that, a linear basis of the polynomial representation
is given in terms of the modified polynomials.
The irreducibility is
 shown by checking that any non-zero vector in the representation is cyclic.
(The notion ``modified polynomials" in this paper
are, in some sense, similar to {\it non-semisimple Macdonald polynomials}
 introduced in \cite{Ch} for reduced root systems.
The feature of the present paper is,
as we stated above,
to obtain an inductive and combinatorial method
by the use of the modified intertwiners.
See Remark \ref{rem:Chered} (ii).)

When the representation is reducible,
we describe subrepresentations
in terms of spanning sets of non-symmetric Koornwinder polynomials,
or we realize them as ideals which are invariant under the action of DAHA.

The main statements in this paper are as follows:

\begin{thm}
(i) The polynomial representation is irreducible and $Y$-semisimple
if the parameters $q,t,a^*,b^*,c^*,d^*$ are not roots of
 the following Laurent polynomials:
\bean
&&t^{(k+1)/m}q^{(r-1)/m}-\omega_m \label{eq:spec_tq_intro}\\
&&\qquad \qquad(n\geq k+1\geq 0, r-1\geq1,m=GCD(k+1,r-1), \nonumber\\
&&\qquad \qquad \mbox{$\omega_m$ is a primitive $m$-th root of unity$)$},
\nonumber
\eean
\bean
&&t^{k+1}q^{r-1}a^*{}^{2}-1 \label{eq:spec_aa_intro}\\
&&\qquad \qquad(2n-2\geq k+1\geq 0, r-1\geq1),\nonumber\\
&&t^{n-i}q^{r-1}a^*b^*{}^{\pm1}-1 \label{eq:spec_ab_intro}\\
&&\qquad \qquad(n\geq i\geq 1, r-1\geq1),\nonumber\\
&&t^{n-i}q^{r-1-\theta(\pm1)}a^*c^*{}^{\pm1}-1 \label{eq:spec_ac_intro}\\
&&\qquad \qquad(n\geq i\geq 1, r-1\geq1),\nonumber\\
&&t^{n-i}q^{r-1-\theta(\pm1)}a^*d^*{}^{\pm1}-1 \label{eq:spec_ad_intro}\\
&&\qquad \qquad(n\geq i\geq 1, r-1\geq1) \nonumber
\eean
where $\theta(+1)=1$ and $\theta(-1)=0$.

(ii) We take a specialization of parameters
 such that the parameters are roots of one of
 (\ref{eq:spec_tq_intro}),
(\ref{eq:spec_aa_intro}),
(\ref{eq:spec_ab_intro}),
(\ref{eq:spec_ac_intro}), or
(\ref{eq:spec_ad_intro})
and are not roots of the others.
Under the specialization, we give a linear basis
 of the polynomial representation
in terms of linear combinations of non-symmetric Koornwinder polynomials.
(The linear combinations are of the form (\ref{eq:mod-poly-form}).
The construction of the basis is given in Theorem \ref{thm:construct-basis}.)

(iii-1) For the case (\ref{eq:spec_tq_intro}) with $k+1\geq2$,
 the polynomial representation is reducible.
It is $Y$-semisimple if $k+1=n$,
and it is not $Y$-semisimple if $2\leq k+1<n$.
We give $N=\lfloor \frac{n}{k+1}\rfloor$ subrepresentations
 $I_1^{(k,r)}\subset I_2^{(k,r)}\subset \cdots \subset I_N^{(k,r)}$
 in terms of vanishing conditions:
 a Laurent polynomial $f(z_1,\ldots,z_n)$
 should be zero if certain ratios or products of $z_i$ take certain values
 (see Definition \ref{defn:wheel_wheel-condition}).
We show that $I_1^{(k,r)}$ is irreducible and $Y$-semisimple.
A basis of $I_1^{(k,r)}$ is given
 (see Theorem \ref{thm:wheel_irr} and Theorem \ref{thm:I1}).
Detailed arguments are given in \S\ref{subsect:wheel}.

(iii-2) For the case (\ref{eq:spec_tq_intro}) with $k+1=1$
 or (\ref{eq:spec_aa_intro}),
 the polynomial representation is irreducible and not $Y$-semisimple.
Detailed arguments are given
 in \S\ref{subsect:tq}, \S\ref{subsect:aa}.

(iii-3) For the case (\ref{eq:spec_ab_intro})
 or (\ref{eq:spec_ac_intro}) or (\ref{eq:spec_ad_intro}),
 the polynomial representation is reducible and $Y$-semisimple.
It has a unique subrepresentation.
If the sign is plus in (\ref{eq:spec_ab_intro}),
 (\ref{eq:spec_ac_intro}), (\ref{eq:spec_ad_intro}),
then the subrepresentation is characterized in terms of a vanishing condition:
a Laurent polynomial $f(z_1,\ldots,z_n)$ should be zero
if $(z_1,\ldots,z_n)$ are included in certain grid points.
This vanishing condition is a generalization of that in \cite{vDSt}.
The explicit statements and detailed arguments are given
 in \S\ref{subsect:ab}, \S\ref{subsect:ac,ad}.

(iii-4) For the case (\ref{eq:spec_tq_intro}) with $k+1=0$
(that is, $q$ is a root of unity),
the polynomial representation is reducible and $Y$-semisimple.
There are infinitely-many subrepresentations labelled by partitions
with length $\leq n$,
and inclusion relation of two subrepresentations
 is determined by the dominance ordering of corresponding partitions.
All the irreducible subquotients
 are finite dimensional, and isomorphic to each other.
Detailed arguments are given in \S\ref{subsect:q}.
\end{thm}

As an application,
 we show some non-symmetric Koornwinder polynomials
 for the specialized parameters
 are eigenfunctions with respect to Demazure-Lusztig operators $\widehat{T}_i$
 for some $i$.
For such Laurent polynomials, one can construct
polynomial solutions of quantum Knizhnik-Zamolodchikov (qKZ)
 equation of type $(C^\vee_n, C_n)$ (about detail, see \cite{KaSh}).
In Proposition 5.4 of \cite{KaTa},
an example of polynomial solutions of qKZ equation (of type $GL_n$)
is given.
In this paper, in Proposition \ref{prop:level1sol},
we give a generalization of the example in \cite{KaTa}.

On the result (iii-1),
there are some variants for symmetric (Laurent) polynomials
 or for the polynomial representation of $GL_n$-DAHA.

In \cite{FJMM}, the terminology ``wheel condition" originally appeared.
They considered a vanishing condition for $\mathfrak{S}_n$-symmetric polynomials
and called it {\it wheel condition}.
They give a linear basis of the ideal defined by the wheel condition
in terms of symmetric Macdonald polynomials for specialized parameters.

In \cite{Ka1}, standing on \cite{FJMM},
the author introduced a vanishing condition
for $BC_n$-symmetric Laurent polynomials,
and give a linear basis of the ideal
in terms of symmetric Macdonald-Koornwinder polynomials
 for specialized parameters.
For $BC_n$-symmetric Laurent polynomials,
the vanishing condition in \cite{Ka1}
is equivalent to $1$-wheel condition in the present paper.

For the DAHA of type $GL_n$,
 the algebra has only two parameters $q$ and $t$.
When the parameters are specialized at $t^{(k+1)/m}q^{(r-1)/m}=\omega_m$
(where $n\geq k+1\geq 2$, $r-1\geq1$, $m=GCD(k+1,r-1)$,
 $\omega_m$ is a primitive $m$-th root of unity),
the author introduced vanishing conditions for Laurent polynomials
 called multi-wheel condition (of type $GL_n$) in \cite{Ka2}.
By this condition,
 a series of subrepresentations of the polynomial representation
 was constructed.
In \cite{En}, it was shown that the series gives the composition series
 of the polynomial representation
 for the case $m=GCD(k+1,r-1)=1$ and $k+1\neq2$.
The result (iii-1) in the present paper
is a $(C^\vee_n, C_n)$-version of the result in \cite{Ka2}.
(The cases (\ref{eq:spec_aa_intro}), $\ldots$, (\ref{eq:spec_ad_intro})
 have no counterparts for $GL_n$-DAHA
 due to the difference of numbers of parameters.)

\medskip

This paper is organized as follows.

In Section 2, we introduce notations in this paper
 and review basic properties for DAHA and its polynomial representation.
Especially, duality of the algebra, the Noumi representation,
 the non-symmetric Koornwinder polynomials,
the intertwiners, irreducibility of the representation,
duality relations and evaluation formulas for the non-symmetric
Koornwinder polynomials are stated.

In Section 3, we consider specializations of parameters
where the representation can be non-$Y$-semisimple or reducible.
We realize such specializations as follows:
First, take an irreducible factor $s$ of
(\ref{eq:spec_tq_intro}),
(\ref{eq:spec_aa_intro}),
(\ref{eq:spec_ab_intro}),
(\ref{eq:spec_ac_intro}), or
(\ref{eq:spec_ad_intro}).
Second, take the quotient of the Laurent polynomial ring of the parameters
by the prime ideal generated by $s$.
Finally, take the fractional field of the quotient ring.
This realization gives a field where the parameters satisfy $s=0$.
For such a specialization, we introduce modified intertwiners,
and in terms of them, we give a method to compute
modified polynomials (generalized $Y$-eigenfunctions).
As a result, a linear basis of the polynomial representation
is given by the modified polynomials.

In Section 4, applying tools given in Section 3,
we examine irreducibility and $Y$-semisimplicity
 of the polynomial representation
for each case of
(\ref{eq:spec_tq_intro}),
(\ref{eq:spec_aa_intro}),
(\ref{eq:spec_ab_intro}),
(\ref{eq:spec_ac_intro}), or
(\ref{eq:spec_ad_intro}).
For the cases where the representation is reducible,
 we describe subrepresentations
 as spanning sets of the non-symmetric Koornwinder polynomials
 or as ideals defined by certain vanishing conditions.
The duality relations and the evaluation formulas
 give values of the non-symmetric Koornwinder polynomial
 at certain grid points.
By using this property,
 we show that some non-symmetric Koornwinder polynomials
satisfy the vanishing conditions.

\medskip

{\bf Acknowledgements.} \quad
The author thanks to his advisor T. Miwa for many comments.
The author also thanks to K. Shigechi for discussions.
The research of the author is partially supported by
Grant-in-Aid for JSPS Fellows (DC1) No.17-02106.


\section{Notations and basic properties}

\subsection{The double affine Hecke algebra}

In this paper, we assume $n\geq2$.
We realize the affine root system of type $C_n$
 in $\mathbb{R}^n\oplus \mathbb{R}\delta$
and denote affine simple roots by:
\bea
&&\alpha_i=\epsilon_i-\epsilon_{i+1} \ (1\leq i\leq n-1),
\quad \alpha_n=2\epsilon_n,
\quad \alpha_0=\delta-2\epsilon_1,
\eea
where $\epsilon_i$ are standard bases:
 $\langle \epsilon_i,\epsilon_j\rangle=\delta_{ij}$
and $\delta$ is a radical element.
Let $\alpha_i^\vee=\frac{2\alpha_i}{\langle \alpha_i,\alpha_i\rangle}$
be simple co-roots and
let $\varpi_i=\epsilon_1+\ldots+\epsilon_i$ be the fundamental weights
and $\wt_n\cong\mathbb{Z}^n$ be the weight lattice.

Let $\mathbb{K}=
\mathbb{C}(q^{1/2},t^{1/2},t_n^{1/2},t_0^{1/2},u_n^{1/2},u_0^{1/2})$.
Denote the ring of $n$-variable Laurent polynomials
by $P_n=\mathbb{K}[x_1^{\pm1},\ldots,x_n^{\pm1}]$.
$P_n$ is isomorphic to the group algebra $\mathbb{K}\wt_n$.
Denote the finite and affine Weyl groups of type $C_n$
by $W_0=\langle s_1,\ldots,s_n \rangle$,
$W=\langle s_0,\ldots,s_n \rangle$, respectively.
The action of the affine Weyl group is given by
$wx^\lambda=x^{w\lambda}$ where $x^\delta=q$.
Especially,
\bea
s_0f(x_1,x_2\ldots)&=&f(qx_1^{-1},x_2,\ldots)\\
s_nf(\ldots,,x_{n-1},x_n)&=&f(\ldots,x_{n-1},x_n^{-1}).
\eea

\begin{defn}[DAHA]\normalfont
The double affine Hecke algebra (DAHA)
 $\HH=\HH(q^{1/2},t^{1/2},t_n^{1/2},t_0^{1/2},u_n^{1/2},u_0^{1/2})$
 of type $(C^\vee_n, C_n)$
 (\cite{Sa}) is
a unital associative $\mathbb{K}$-algebra generated by
 $T_0,\ldots,T_{n}, X_1^{\pm1},\ldots,X_n^{\pm1}$
 with defining relations as follows:
\bea
\mbox{quadratic Hecke relations:}&& \\
(T_0-t_0^{1/2})(T_0+t_0^{-1/2})&=&0, \\
(T_i-t^{1/2})(T_i+t^{-1/2})&=&0 \qquad 1\leq i\leq n-1, \\
(T_n-t_n^{1/2})(T_n+t_n^{-1/2})&=&0,
\eea
\bea
\mbox{braid relations:} \\
T_0T_1T_0T_1&=&T_1T_0T_1T_0, \\
T_iT_{i+1}T_i&=&T_{i+1}T_iT_{i+1} \qquad 1\leq i\leq n-2, \\
T_{n-1}T_nT_{n-1}T_n&=&T_nT_{n-1}T_nT_{n-1}, \\
T_iT_j&=&T_jT_i \qquad |i-j|\geq2,
\eea
\bea
\mbox{relations between $X$ and $T$:} \\
X_iX_j&=&X_jX_i \quad \forall i,j, \\
T_iX_j&=&X_jT_i \quad\langle\alpha_i,\epsilon_j\rangle=0, \\
T_iX_iT_i&=&X_{i+1} \quad{1\leq i\leq n-1}, \\
X_n^{-1}T_n^{-1}&=&T_nX_n+(u_n^{1/2}-u_n^{-1/2}), \\
q^{-1/2}T_0^{-1}X_1&=&q^{1/2}X_1^{-1}T_0+(u_0^{1/2}-u_0^{-1/2}).
\eea
\end{defn}

The subalgebra
 $\HH^{\mathrm{aff}}=\HH^{\mathrm{aff}}(t^{1/2},t_n^{1/2}t_0^{1/2})\subset\HH$
 generated by $T_0,\ldots,T_n$ is
the affine Hecke algebra of type $C_n$.
Let $\HH^{Y}=\HH^{Y}(t^{1/2},t_n^{1/2}t_0^{1/2})$
 be the algebra generated by $T_1,\ldots,T_n$ and $Y_1^{\pm1},\ldots,Y_n^{\pm1}$
with the defining relations as follows:
\bea
\mbox{the quadratic Hecke relations}&\mbox{for}&\mbox{$T_1,\ldots,T_n$}, \\
\mbox{the braid relations}&\mbox{for}&\mbox{$T_1,\ldots,T_n$}, \\
Y_iY_j&=&Y_jY_i \quad \forall i,j, \\
T_iY_j&=&Y_jT_i \quad\langle\alpha_i,\epsilon_j\rangle=0, \\
T_iY_{i+1}T_i&=&Y_i \quad{1\leq i\leq n-1}, \\
T_n^{-1}Y_n&=&Y_n^{-1}T_n+(t_0^{1/2}-t_0^{-1/2}).
\eea
Then $\HH^{Y}$ is isomorphic to $\HH^{\mathrm{aff}}$
by the correspondence
\bea
Y_i &\mapsto& T_i \ldots T_{n-1}T_n \ldots T_0 T_1^{-1} \ldots T_{i-1}^{-1}.
\eea
Hereafter we identify $Y_i$ with the right hand side above.
We see that the elements $T_i$ and $Y_j$ satisfy the following relation:
\bea
q^{-1}Y_1^{-1}U_n^{-1}&=&U_nY_1+q^{-1/2}(u_0^{1/2}-u_0^{-1/2}) \\
\mbox{where}&& \\
U_n&:=&X_1^{-1}T_0Y_1^{-1}.
\eea

Let $*:\mathbb{K}\rightarrow\mathbb{K}$ be an involution
given by
\bea
t_0^{1/2}\leftrightarrow u_n^{1/2}
\eea
and the other parameters $t^{1/2},q^{1/2},t_n^{1/2},u_0^{1/2}$ are fixed.
There is an anti-automorphism of DAHA:
\bea
*:\HH(q^{1/2},t^{1/2},t_n^{1/2},t_0^{1/2},u_n^{1/2},u_0^{1/2})\rightarrow
\HH(q^{1/2},t^{1/2},t_n^{1/2},u_n^{1/2},t_0^{1/2},u_0^{1/2})
\eea
given by
\bea
T_0\mapsto U_n,\ T_i\mapsto T_i,\ X_i\mapsto Y_i^{-1},\ Y_i\mapsto X_i^{-1}
\quad (i=1,\ldots,n).
\eea
The anti-automorphism $*$ is called duality anti-involution.
Especially, the subalgebra $\HH^{X}\subset\HH$ generated by
$T_1,\ldots,T_n$ and $X_1,\ldots,X_n$ is anti-isomorphic to
the algebra $\HH^{Y}(t^{1/2},t_n^{1/2}u_n^{1/2})$.
The parameters
\bea
(q^{1/2}{}^*,t^{1/2}{}^*,t_n^{1/2}{}^*,t_0^{1/2}{}^*,u_n^{1/2}{}^*,u_0^{1/2}{}^*)
=(q^{1/2},t^{1/2},t_n^{1/2},u_n^{1/2},t_0^{1/2},u_0^{1/2})
\eea
are called dual parameters.

We often use another notations for parameters (\cite{Sa}):
\bea
&&a=t_n^{1/2}u_n^{1/2},
b=-t_n^{1/2}u_n^{-1/2},
c=q^{1/2}t_0^{1/2}u_0^{1/2},
d=-q^{1/2}t_0^{1/2}u_0^{-1/2}. 
\eea
\bea
&&a'=t_n^{-1/2}t_0^{-1/2},
b'=-t_n^{-1/2}t_0^{1/2},
c'=q^{-1/2}u_n^{-1/2}u_0^{-1/2},
d'=-q^{-1/2}u_n^{-1/2}u_0^{1/2}.
\eea
\bea
&&a^*=t_n^{1/2}t_0^{1/2},
b^*=-t_n^{1/2}t_0^{-1/2},
c^*=q^{1/2}u_n^{1/2}u_0^{1/2},
d^*=-q^{1/2}u_n^{1/2}u_0^{-1/2}. 
\eea
Note that $a^*,\ldots,d^*$ are dual parameters of $a,\ldots,d$,
and $a',\ldots,d'$ are inverses of $a^*,\ldots,d^*$.

An action of $\HH$
is realized on the ring of $n$-variable Laurent polynomials $P_n$.
\begin{prop}[polynomial representation]
Define the linear operators on $P_n$ (\cite{No}):
\bea
\widehat{T}_0^{\pm1}&=&t_0^{\pm1/2}+t_0^{-1/2}
\frac{(1-cx_1^{-1})(1-dx_1^{-1})}{1-qx_1^{-2}}(s_0-1)\\
\widehat{T}_i^{\pm1}&=&t^{\pm1/2}+t^{-1/2}
\frac{1-tx_ix_{i+1}^{-1}}{1-x_ix_{i+1}^{-1}}(s_i-1)\\
\widehat{T}_n^{\pm1}&=&t_n^{\pm1/2}+t_n^{-1/2}
\frac{(1-ax_n)(1-bx_n)}{1-x_n^2}(s_n-1).
\eea
The map $T_i\mapsto \widehat{T}_i$ and $X_j\mapsto x_j$
$(0\leq i\leq n,1\leq j \leq n)$
gives a representation of DAHA (\cite{Sa}).
This is called the polynomial representation, or the Noumi representation.
\end{prop}

\subsection{The non-symmetric Koornwinder polynomials}

It is known that
the q-Dunkl-Cherednik operators $\widehat{Y}_i$ given by
\bea
\widehat{Y}_i:=\widehat{T}_i \ldots \widehat{T}_{n-1}\widehat{T}_n
 \ldots \widehat{T}_0 \widehat{T}_1^{-1} \ldots \widehat{T}_{i-1}^{-1}
\eea
 are mutually commutative and
 triangular with respect to an partial ordering
and the pairs of eigenvalues of $\widehat{Y}_1,\ldots,\widehat{Y}_n$
 are mutually different.
Consequently,
 the space $P_n$ is simultaneously diagonalizable
 with respect to $\widehat{Y}_1,\ldots,\widehat{Y}_n$.
Their joint eigenfunctions are called non-symmetric Koornwinder polynomials.
We call this property $Y$-semisimplicity.

We will explain the explicit definition of the ordering, and the eigenvalues,
and so on.

Denote by $\lambda^+$ the unique dominant element in $W_0\lambda$.
($\lambda^+$ is a partition of length $\leq n$.)
Take the shortest element $w\in W_0$ such that
 $w\lambda^+=\lambda$ and denote it by $w_\lambda^+$.
Put $\rho=(n-1,n-2,\ldots,1,0)$, $\rho(\lambda)=w_\lambda^+\rho$,
$\sigma(\lambda)=
\left(\sgn(\lambda_1),\ldots,\sgn(\lambda_n)\right)$
where $\sgn(0)=+1$.
We identify any element $w\in W_0$
with a permutation $\{\pm1,\ldots,\pm n\}\rightarrow\{\pm1,\ldots,\pm n\}$
by $\pm i\mapsto\langle (\sum_{j=1}^n  j\epsilon_j),\pm w\epsilon_i\rangle$.
Then we see that $\sgn(w_\lambda^+(i))=\sgn(\lambda_i)$.
We extend the action of $W_0$ on $\mathbb{R}^n$
 to an action of $W$ on $\mathbb{R}^n$ by
\bea
s_0\cdot(v_1,v_2,\ldots,v_n)=(-1-v_1,v_2,\ldots,v_n).
\eea
We call it the {\it dot} action of $W$.

Define partial orderings $\lambda\geq\mu$ and $\lambda\succeq\mu$
 in $\mathbb{Z}^n$ as follows:
\bea
\lambda\geq\mu \quad&\mbox{if}&
\mbox{$\lambda-\mu\in\sum_{i=1}^n \mathbb{Z}_{\geq0}\alpha_i^\vee$}, \\
\lambda\succeq\mu \quad&\mbox{if}&
\mbox{$\lambda^+>\mu^+$, or $\lambda^+=\mu^+$ and $\lambda\geq\mu$}.
\eea

\begin{defn}[non-symmetric Koornwinder polynomials (\cite{No,Sa,St})]\normalfont
For $\lambda\in\mathbb{Z}^n$, the non-symmetric Koornwinder polynomial
$E_\lambda$ is defined by
\bea
\widehat{Y}_iE_\lambda
 &=&y(\lambda)_iE_\lambda \\
E_\lambda&=&x^\lambda+\sum_{\mu\prec\lambda} c_{\lambda\mu}x^\mu \qquad (c_{\lambda\mu}\in\mathbb{K})
 \\
\mbox{where}&&\\
\quad y(\lambda)_i&:=&q^{\lambda_i}t^{\rho(\lambda)_i}a^*{}^{\sigma(\lambda)_i}.
\eea
\end{defn}
For the pair of eigenvalues $y(\lambda)=(y(\lambda)_1,\ldots,y(\lambda)_n)$,
and any Laurent polynomial (or rational function) $f$,
we express the scalar $f(y(\lambda))\in\mathbb{K}$
 by $f(Y)|_\lambda$ or $f(\lambda)$:
\bea
f(Y)|_\lambda &\stackrel{\mathrm{def}}{=}& f(y(\lambda)) \\
f(\lambda) &\stackrel{\mathrm{def}}{=}& f(y(\lambda)).
\eea
We denote by $E_\lambda^*$ the polynomial
 where the parameters of $E_\lambda$ are replaced by
 the dual parameters.
We see that $\{E_\lambda;\lambda\in\mathbb{Z}^n\}$
or $\{E_\lambda^*;\lambda\in\mathbb{Z}^n\}$ constitutes
a $\mathbb{K}$-basis of $P_n$.

\subsection{The intertwiners}

There are operators which send an $Y$-eigenfunction to other $Y$-eigenfunctions.
Such operators are called intertwining operators or intertwiners.

\begin{defn}[intertwiners]\label{defn:int}
Let $\phi_i$ be as follows:
\bea
\phi_i&=&T_i+\frac{t^{1/2}-t^{-1/2}}{Y_{i+1}/Y_i-1} \quad(1\leq i\leq n-1), \\
\phi_n&=&T_n+\frac{(t_n^{1/2}-t_n^{-1/2})+(t_0^{1/2}-t_0^{-1/2})Y_n^{-1}}
{Y_n^{-2}-1}, \\
\phi_0&=&U_n+\frac{(u_n^{1/2}-u_n^{-1/2})+(u_0^{1/2}-u_0^{-1/2})q^{1/2}Y_1}
{qY_1^2-1}.
\eea
\end{defn}
Since there are rational functions $f(Y)$ in the right hand sides,
these operators $\phi_i$ are not elements of $\HH$.
However $\phi_i$ are considered as linear operators in $P_n$
by replacing $f(Y)$ by $f(Y)|_\lambda$ on each eigenspace $\mathbb{K}E_\lambda$.

For any $\lambda+m\delta\in\mathbb{Z}^n\oplus\frac{1}{2}\mathbb{Z}\delta$,
we define
$Y^{\lambda+m\delta}=q^{-m}Y_1^{\lambda_1}\cdots Y_n^{\lambda_n}$.
Define $1$-variable Laurent polynomials $D_i$ and $N_i$ ($1\leq i\leq n$) by
$D_i(x)=x^{-1}-1$ ($0\leq i\leq n$) and
\bea
N_i(x)
&=&t^{1/2}(x^{-1}-t^{-1}) \qquad(1\leq i\leq n-1), \\
N_n(x)
&=&t_n^{1/2}(x^{-1}-a')(x^{-1}-b'), \\
N_0(x)
&=&u_n^{1/2}(x^{-1}-q^{1/2}c')(x^{-1}-q^{1/2}d').
\eea
Then $D_i(Y^{\pm\alpha_i^\vee})$ and $N_i(Y^{\pm\alpha_i^\vee})$
 are given as follows:
\bean
D_i(Y^{\pm\alpha_i^\vee})&=&Y_{i+1}^{\pm1}Y_i^{\mp1}-1
 \qquad(1\leq i\leq n-1) \nonumber\\
D_n(Y^{\pm\alpha_n^\vee})&=&Y_n^{\mp2}-1 \label{eq:def-of-D}\\
D_0(Y^{\pm\alpha_0^\vee})&=&q^{\pm}Y_1^{\pm2}-1\nonumber
\eean
and
\bean
N_i(Y^{\pm\alpha_i^\vee})
&=&t^{1/2}(Y_{i+1}^{\pm1}Y_i^{\mp1}-t^{-1})
 \qquad(1\leq i\leq n-1) \nonumber\\
N_n(Y^{\pm\alpha_n^\vee})
&=&t_n^{1/2}(Y_n^{\mp1}-a')(Y_n^{\mp1}-b') \label{eq:def-of-N}\\
N_0(Y^{\pm\alpha_0^\vee})
&=&u_n^{1/2}(q^{\pm1/2}Y_1^{\pm1}-q^{1/2}c')(q^{\pm1/2}Y_1^{\pm1}-q^{1/2}d').
\nonumber
\eean
Note that $D_i(Y^{\alpha_i^\vee})$ is the denominator of $\phi_i$.

\begin{prop}\label{prop:int}
For any $0\leq i\leq n$, we have
\bea
\phi_iY^{\epsilon_j}&=&Y^{s_i\epsilon_j}\phi_i \quad (1\leq j\leq n),\\
\phi_i^2&=&\frac{N_i(Y^{\alpha_i^\vee})N_i(Y^{-\alpha_i^\vee})}
{D_i(Y^{\alpha_i^\vee})D_i(Y^{-\alpha_i^\vee})}, \\
\phi_{i}\phi_{i+1}\phi_{i}&=&\phi_{i+1}\phi_{i}\phi_{i+1}
 \quad(1\leq i\leq n-2), \\
\phi_{i}\phi_{i+1}\phi_{i}\phi_{i+1}&=&
\phi_{i+1}\phi_{i}\phi_{i+1}\phi_{i} \quad(i=0,n-1).
\eea
We see that $\phi_iE_\lambda$ is proportional to $E_{s_i\cdot\lambda}$.
Let $c_{i,\lambda}$ be the coefficient
\bea
\phi_iE_\lambda=c_{i,\lambda}E_{s_i\cdot\lambda}.
\eea
Then for $1\leq i \leq n-1$,
\bea
c_{i,\lambda}\ =&t^{1/2}&
 \qquad\mbox{if $\langle\lambda,\alpha_i\rangle<0$}, \\
c_{i,\lambda}\ =&0&
 \qquad\mbox{if $\langle\lambda,\alpha_i\rangle=0$}, \\
c_{i,\lambda}\ =&t^{-1/2}\phi_i^2|_\lambda&
 \qquad\mbox{if $\langle\lambda,\alpha_i\rangle>0$},
\eea
and
\bea
c_{n,\lambda}\ =&t_n^{1/2}&
 \qquad\mbox{if $\langle\lambda,\alpha_n\rangle<0$}, \\
c_{n,\lambda}\ =&0&
 \qquad\mbox{if $\langle\lambda,\alpha_n\rangle=0$}, \\
c_{n,\lambda}\ =&t_n^{-1/2}\phi_n^2|_\lambda&
 \qquad\mbox{if $\langle\lambda,\alpha_n\rangle>0$},
\eea
\bea
c_{0,\lambda}\ =&t_0^{1/2}t^{-\rho(\lambda)_1}a^{*}{}^{-1}&
 \qquad\mbox{if $\langle\lambda,\alpha_0\rangle\leq0$, that is $\lambda_1\geq0$},\\
c_{0,\lambda}\ =&t_0^{-1/2}t^{\rho(\lambda)_1}a^{*}\phi_0^2|_\lambda&
 \qquad\mbox{if $\langle\lambda,\alpha_0\rangle>0$, that is $\lambda_1<0$}.
\eea
\end{prop}
Note that $c_{i,\lambda}$ are Laurent monomials in $\mathbb{K}$
if $\langle\lambda,\alpha_i\rangle\leq0$
and $c_{i,\lambda}$ are rational functions in $\mathbb{K}$
if $\langle\lambda,\alpha_i\rangle>0$.

\begin{rem}\normalfont
Originally, Sahi introduced the intertwining operators $S_i\in\HH$ 
for type $(C^\vee_n,C_n)$ (\cite{Sa}) by
\bea
S_i=[T_i,Y_i],\ S_n=[T_n,Y_n],\ S_0=[Y_1,U_n].
\eea
The relations between $S_i$ and $\phi_i$ are given by
$S_i=\phi_i\cdot(Y_i-Y_{i+1})$, $S_n=\phi_n\cdot(Y_n-Y_n^{-1})$,
 $S_0=(Y_1-q^{-1}Y_1^{-1})\cdot\phi_0$.
In \cite{St} or \cite{NUKW}, similar intertwiners $S_i$ are also introduced.
However, the coefficient $c_{i,\lambda}'$
in $S_iE_\lambda=c_{i,\lambda}'E_{s_i\cdot\lambda}$
becomes more complicated even if $\langle\lambda,\alpha_i\rangle<0$.
We introduce $\phi_i$ in order to take the coefficient
$c_{i,\lambda}$ as a Laurent monomial in $\mathbb{K}$
 if $\langle\lambda,\alpha_i\rangle<0$.
\end{rem}

Since the dot action of $W$ on $\mathbb{Z}^n$ is transitive,
 we can check that any $E_\lambda$ ($\lambda\in\mathbb{Z}^n$)
 is a cyclic vector in $P_n$
by applying the intertwiners $\phi_i$.
Hence, we obtain the statement as follows:
\begin{prop}[irreducibility(\cite{Sa})]
The polynomial representation $P_n$ is irreducible.
\end{prop}

\subsection{The duality relations and the evaluation formulas}

We introduce two properties called duality relations and evaluation formulas
for $E_\lambda$.

\begin{defn}\normalfont
For $\lambda\in\mathbb{Z}^n$ and $f\in P_n$,
 we define two maps $\chi_\lambda$ and $\chi_\lambda^*$:
 $P_n\rightarrow \mathbb{K}$ as follows:
\bea
\chi_\lambda(f)&=& 
f(y(\lambda)_1^{-1},\cdots,y(\lambda)_n^{-1}) \\
\chi_\lambda^*(f)&=& 
f(y^*(\lambda)_1^{-1},\cdots,y^*(\lambda)_n^{-1}).
\eea
\end{defn}

\begin{prop}[duality relations(\cite{Sa})]\label{prop:duality}
For any $\lambda,\mu\in\mathbb{Z}^n$, we have
\bea
\chi_\mu^*(E_\lambda)\chi_0(E_\mu^*)
=\chi_\lambda(E_\mu^*)\chi_0^*(E_\lambda).
\eea
\end{prop}

We have relations between
 $\chi_0^*(E_{s_i\cdot\lambda})$ and $\chi_0^*(E_{\lambda})$.

\begin{lem}[recurrence relations]\label{lem:chi0_recur}
(Recall (\ref{eq:def-of-D}) and (\ref{eq:def-of-N})
for the definition of $D_i$ and $N_i$.)
If $\langle\lambda,\alpha_i\rangle<0$ $(1\leq i\leq n-1)$, then
\bea
\chi_0^*(E_{s_i\cdot\lambda})=t^{-1/2}\left.\frac{N_i(Y^{\alpha_i^\vee})}
{D_i(Y^{\alpha_i^\vee})}\right|_{\lambda} \chi_0^*(E_\lambda).
\eea
If $\langle\lambda,\alpha_n\rangle<0$, then
\bea
\chi_0^*(E_{s_n\cdot\lambda})
&=&t_n^{-1/2}\left.\frac{N_n(Y^{\alpha_n^\vee})}
{D_n(Y^{\alpha_n^\vee})}\right|_{\lambda} \chi_0^*(E_\lambda).
\eea
If $\langle\lambda,\alpha_0\rangle<0$, (that is $\lambda_1\geq 0$,) then
\bea
\chi_0^*(E_{s_0\cdot\lambda})
&=&t_0^{-1/2}t^{\rho(\lambda)_1}a^{*}
\left.\frac{N_0(Y^{\alpha_0^\vee})}
{D_0(Y^{\alpha_0^\vee})}\right|_{\lambda} \chi_0^*(E_\lambda).
\eea
\end{lem}
\begin{proof}
The statements are shown by computing
 $\chi_0^*(\phi_iE_\lambda)$ ($0\leq i \leq n$).
\end{proof}

\begin{prop}[evaluation formula(\cite{St})]\label{prop:chi0_formula}
Suppose $\lambda=\lambda^+$. Then we have
\bea
\chi_0^*(E_\lambda)&=&
\prod_{i<j}t^{-(\lambda_i-\lambda_j)}
\frac{(t^{j-i+1}q;q)_{\lambda_i-\lambda_j}}
{(t^{j-i}q;q)_{\lambda_i-\lambda_j}} \\
&&\times\prod_{i<j}t^{-(\lambda_i+\lambda_j)}
\frac{(t^{2n-i-j+1}a^*{}^2q;q)_{\lambda_i+\lambda_j}}
{(t^{2n-i-j}a^*{}^2q;q)_{\lambda_i+\lambda_j}} \\
&&\times\prod_{i}a^{-\lambda_i}t^{(n-i)\lambda_i}
\frac{(t^{n-i}a^*{}^2q,t^{n-i}a^*b^*q,t^{n-i}a^*c^*,t^{n-i}a^*d^*;q)_{\lambda_i}}
{(q^2t^{2(n-i)}a^*{}^2,qt^{2(n-i)}a^*{}^2;q^2)_{\lambda_i}},
\eea
where $(x;q)_i=\prod_{l=0}^{i-1}(1-xq^l)$
 and $(x_1,\ldots,x_j;q)_i=\prod_{l=1}^j(x_l;q)_i$.
\end{prop}
\begin{proof}
This formula is true for $\lambda=(0,\ldots,0)$.
We can show that it is true for any $\lambda=\lambda^+$ inductively,
 using the recurrence relations for $\chi_0^*(E_\lambda)$.
\end{proof}


\section{Specialization of parameters}\label{sect:spec}
In this section, we introduce a specialization of parameters
where the polynomial representation can be non-$Y$-semisimple or reducible.
Put $\mathcal{A}=\mathbb{C}[q^{\pm1/2},t^{\pm1/2},t_n^{\pm1/2},t_0^{\pm1/2},
u_n^{\pm1/2},u_0^{\pm1/2}]$.
($\mathbb{K}$ is the fractional field of $\mathcal{A}$.)

\subsection{Generic and specialized parameters}
We consider a specialization of parameters
 by replacing the field $\mathbb{K}$ by some field including $\mathbb{C}$.
Throughout this paper, we always assume that
$q^{1/2},t^{1/2},t_n^{1/2},t_0^{1/2},
u_n^{1/2},u_0^{1/2}\neq0$
under any specialization of the parameters.

\begin{prop}[generic parameters]\label{prop:generic-param}
The polynomial representation for specialized parameters
is $Y$-semisimple and irreducible
 if the parameters are not roots of either of Laurent polynomials
 given as follows:
\bean
&&t^{(k+1)/m}q^{(r-1)/m}-\omega_m \label{eq:spec_tq}\\
&&\qquad \qquad(n\geq k+1\geq 0, r-1\geq1,m=GCD(k+1,r-1), \nonumber\\
&&\qquad \qquad \mbox{$\omega_m$ is a primitive $m$-th root of unity$)$},
\nonumber\\
&&t^{k+1}q^{r-1}a^*{}^{2}-1 \label{eq:spec_aa}\\
&&\qquad \qquad(2n-2\geq k+1\geq 0, r-1\geq1),\nonumber\\
&&t^{n-i}q^{r-1}a^*b^*{}^{\pm1}-1 \label{eq:spec_ab}\\
&&\qquad \qquad(n\geq i\geq 1, r-1\geq1),\nonumber\\
&&t^{n-i}q^{r-1-\theta(\pm1)}a^*c^*{}^{\pm1}-1 \label{eq:spec_ac}\\
&&\qquad \qquad(n\geq i\geq 1, r-1\geq1),\nonumber\\
&&t^{n-i}q^{r-1-\theta(\pm1)}a^*d^*{}^{\pm1}-1 \label{eq:spec_ad}\\
&&\qquad \qquad(n\geq i\geq 1, r-1\geq1). \nonumber
\eean
\end{prop}

Note that $a^*{}^2=t_nt_0$, $a^*b^*=-t_n$, $a^*b^*{}^{-1}=-t_0$, and
\bea
q^{-\theta(\pm1)}a^*c^*{}^{\pm1}
&=&q^{-1/2}t_n^{1/2}t_0^{1/2}u_n^{\pm1/2}u_0^{\pm1/2}, \\
q^{-\theta(\pm1)}a^*d^*{}^{\pm1}
&=&-q^{-1/2}t_n^{1/2}t_0^{1/2}u_n^{\pm1/2}u_0^{\mp1/2}.
\eea

In order to prove Proposition \ref{prop:generic-param},
we give a sufficient condition where all $Y$-eigenvalues do not degenerate.

\begin{lem}\label{lem:y-coincide}
Suppose that $y(\lambda)=y(\mu)$ under a specialization of parameters
for some $\lambda\neq\mu\in\mathbb{Z}^n$.
Then the parameters are roots of either of the following Laurent polynomials:
\bean
&&t^{k+1}q^{r-1}-1\ \quad(n-1\geq k+1\geq 0, r-1\geq1),
\ \label{eq:spec_zero-of-D_tq}\\
&&t^{k+1}q^{r-1}a^*{}^2-1\  \quad(2n-2\geq k+1\geq 0, r-1\geq1).\qquad
\label{eq:spec_zero-of-D_aa}
\eean
\end{lem}
\begin{proof}
We introduce a notation 
${\bf y}(\lambda)_i:=(\lambda_i,\rho(\lambda)_i,\sigma(\lambda)_i)$
for any $\lambda\in\mathbb{Z}^n$.
Suppose that $y(\lambda)=y(\mu)$ under a specialization of parameters
for some $\lambda\neq\mu$.
Let $i_1,\ldots,i_n$ be
\bea
(i_1,\ldots,i_n)=(|w_\lambda^+(1)|,\ldots,|w_\lambda^+(n)|).
\eea
Then, $|\rho(\lambda)_{i_m}|=n-i_m$ for any $1\leq m\leq n$.
Fix $1\leq\ell\leq n$ such that
${\bf y}(\lambda)_{i_\ell}\neq{\bf y}(\mu)_{i_\ell}$
and ${\bf y}(\lambda)_{i_m}={\bf y}(\mu)_{i_m}$ for any $1\leq m< \ell$.

We divide the proof by the signs of
$(\sigma(\lambda)_{i_\ell},\sigma(\mu)_{i_\ell})$.

(i) If $(\sigma(\lambda)_{i_\ell},\sigma(\mu)_{i_\ell})=(+,-)$,
then the condition $y(\lambda)_{i_\ell}=y(\mu)_{i_\ell}$ implies that
\bea
q^{\lambda_{i_\ell}-\mu_{i_\ell}}
t^{\rho(\lambda)_{i_\ell}-\rho(\mu)_{i_\ell}}a^*{}^2=1.
\eea
Since $(\lambda_{i_\ell})-(\mu_{i_\ell})\geq 0-(-1)=1$ and
$2n-2\geq\rho(\lambda)_{i_\ell}-\rho(\mu)_{i_\ell}\geq0$,
the parameters are roots of (\ref{eq:spec_zero-of-D_aa}).

(ii) The proof for the case
 $(\sigma(\lambda)_{i_\ell},\sigma(\mu)_{i_\ell})=(-,+)$
is similar to that for the case $(+,-)$.

(iii) 
If $(\sigma(\lambda)_{i_\ell},\sigma(\mu)_{i_\ell})=(+,+)$,
then the condition $y(\lambda)_{i_\ell}=y(\mu)_{i_\ell}$ implies that
\bean
q^{\lambda_{i_\ell}-\mu_{i_\ell}}
t^{\rho(\lambda)_{i_\ell}-\rho(\mu)_{i_\ell}}=1. \label{eq:y-coincide_proof1}
\eean
By the definition of $i_\ell$,
we see $n-1\geq\rho(\lambda)_{i_\ell}-\rho(\mu)_{i_\ell}\geq0$.
If $\rho(\lambda)_{i_\ell}-\rho(\mu)_{i_\ell}=0$, then
since ${\bf y}(\lambda)_{i_\ell}\neq{\bf y}(\mu)_{i_\ell}$,
we see that $\lambda_{i_\ell}\neq\mu_{i_\ell}$.
Thus (\ref{eq:y-coincide_proof1}) implies that
the parameters are roots of (\ref{eq:spec_zero-of-D_tq}).
Assume that $n-1\geq\rho(\lambda)_{i_\ell}-\rho(\mu)_{i_\ell}\geq1$.
If $\lambda_{i_\ell}-\mu_{i_\ell}\geq1$,
 then (\ref{eq:y-coincide_proof1}) implies that
 the parameters are roots of (\ref{eq:spec_zero-of-D_tq}).

Hence assume that $\lambda_{i_\ell}\leq\mu_{i_\ell}$.
Take the index $j$ such that $|\rho(\mu)_j|=n-i_\ell$.
Since $|\rho(\mu)_j|=n-i_\ell=\rho(\lambda)_{i_\ell}>\rho(\mu)_{i_\ell}$,
 we see that $|\mu_j|\geq \mu_{i_\ell}$ and $j\neq i_\ell$.
Since $|\rho(\lambda)_j|< n-i_\ell=\rho(\lambda)_{i_\ell}$,
 we have $|\lambda_j|\leq \lambda_{i_\ell}$.
Therefore, $|\lambda_j|\leq \lambda_{i_\ell}\leq\mu_{i_\ell}\leq |\mu_j|$.
If the sign $\sigma(\lambda)_j$ is opposite to the sign $\sigma(\mu)_j$,
then by the same argument above,
the condition $y(\lambda)_j=y(\mu)_j$ implies that
the parameters are roots of (\ref{eq:spec_zero-of-D_aa}).
If $(\sigma(\lambda)_j,\sigma(\mu)_j)=(+,+)$,
then
the condition $y(\lambda)_j=y(\mu)_j$ implies that
\bea
q^{\mu_j-\lambda_j}t^{\rho(\mu)_j-\rho(\lambda)_j}=1
\eea
where $\mu_j-\lambda_j>0$ and $n-1\geq\rho(\mu)_j-\rho(\lambda)_j>0$.
Therefore the parameters are roots of (\ref{eq:spec_zero-of-D_tq}).
If $(\sigma(\lambda)_j,\sigma(\mu)_j)=(-,-)$, 
then the condition $y(\lambda)_j=y(\mu)_j$ implies that
\bea
q^{\lambda_j-\mu_j}t^{\rho(\lambda)_j-\rho(\mu)_j}=1
\eea
where $\lambda_j-\mu_j>0$ and $n-1\geq\rho(\lambda)_j-\rho(\mu)_j>0$.
Therefore the parameters are roots of (\ref{eq:spec_zero-of-D_tq}).

(iv)
The proof for the case $(\sigma(\lambda)_{i_\ell},\sigma(\mu)_{i_\ell})=(-,-)$
 is similar to that for the case $(+,+)$.
\end{proof}

We will start the proof of Proposition \ref{prop:generic-param}.
For parameters satisfying the assumption,
the $Y$-semisimplicity and irreducibility of $P_n$ are shown by the following
statements: for any $\lambda\in\mathbb{Z}^n$,
\begin{itemize}
\item
$y(\lambda)\neq y(\mu)$ for any $\mu\in\mathbb{Z}^n$ such that $\mu\neq\lambda$
and $E_\lambda$ is well-defined under the specialization of parameters.
\item
$E_{s_i\cdot\lambda}\in \HH E_\lambda$.
\item
The dot-action $\cdot$ of $W$ on $\mathbb{Z}^n$ is transitive.
\end{itemize}
The first statement is already shown by Lemma \ref{lem:y-coincide},
and the third statement is easy.

\begin{proof}[Proof of Proposition \ref{prop:generic-param}]

(See (\ref{eq:def-of-D}) and (\ref{eq:def-of-N}) for the definition of
 $D_i(Y^{\pm\alpha_i^\vee})$ and $N_i(Y^{\pm\alpha_i^\vee})$.)

For any $0\leq i \leq n$ and any $\lambda\in\mathbb{Z}^n$
 such that $\lambda\neq s_i\cdot\lambda$,
each factor in $D_i(Y^{\alpha_i^\vee})D_i(Y^{-\alpha_i^\vee})|_\lambda$
 is either (\ref{eq:spec_zero-of-D_tq}) or (\ref{eq:spec_zero-of-D_aa})
up to a monic Laurent monomial.
We see that
(\ref{eq:spec_zero-of-D_tq}) or (\ref{eq:spec_zero-of-D_aa})
is factorized in $\mathcal{A}$ into (\ref{eq:spec_tq}), (\ref{eq:spec_aa}),
 respectively.

For any $0\leq i \leq n$ and any $\lambda\in\mathbb{Z}^n$
 such that $\lambda\neq s_i\cdot\lambda$,
each factor in $N_i(Y^{\alpha_i^\vee})N_i(Y^{-\alpha_i^\vee})|_\lambda$
 is either of the following Laurent polynomials
up to a monic Laurent monomial:
\bean
t^{k+1}q^{r-1}-1&& \quad(n\geq k+1\geq 0, r-1\geq1),
\ \mbox{or}\label{eq:spec_zero-of-N_tq}\\
t^{k+1}q^{r-1}a^*{}^2-1&& \quad(2n-2\geq k+1\geq 0, r-1\geq1),
\ \mbox{or}\label{eq:spec_zero-of-N_aa}\\
t^{n-i}q^{r-1}a^*b^*{}^{\pm1}-1&& \quad(n\geq i\geq 1, r-1\geq1),
\ \mbox{or}\label{eq:spec_zero-of-N_ab}\\
t^{n-i}q^{r-1-\theta(\pm1)}a^*c^*{}^{\pm1}-1&& \quad(n\geq i\geq 1, r-1\geq1),
\ \mbox{or}\label{eq:spec_zero-of-N_ac}\\
t^{n-i}q^{r-1-\theta(\pm1)}a^*d^*{}^{\pm1}-1&& \quad(n\geq i\geq 1, r-1\geq1),
\label{eq:spec_zero-of-N_ad}
\eean
where $\theta(+1)=1$ and $\theta(-1)=0$.
We see that
(\ref{eq:spec_zero-of-N_tq}),
(\ref{eq:spec_zero-of-N_aa}),
(\ref{eq:spec_zero-of-N_ab}),
(\ref{eq:spec_zero-of-N_ac}), or
(\ref{eq:spec_zero-of-N_ad})
 is factorized in $\mathcal{A}$ into (\ref{eq:spec_tq}), (\ref{eq:spec_aa}),
(\ref{eq:spec_ab}), (\ref{eq:spec_ac}), or (\ref{eq:spec_ad}), respectively.

Suppose that specialized parameters are
 not roots of either of (\ref{eq:spec_tq}),
$\ldots$, (\ref{eq:spec_ad}).
Then we have
\begin{enumerate}
\item $N_i(Y^{\alpha_i^\vee})N_i(Y^{-\alpha_i^\vee})|_\lambda\neq0$
for any $0\leq i\leq n$ and $\lambda\in\mathbb{Z}^n$ such that
$\lambda\neq s_i\cdot\lambda$,
\item $D_i(Y^{\alpha_i^\vee})D_i(Y^{-\alpha_i^\vee})|_\lambda\neq0$
for any $0\leq i\leq n$ and $\lambda\in\mathbb{Z}^n$ such that
$\lambda\neq s_i\cdot\lambda$,
\item the specialized parameters are not any roots of either
(\ref{eq:spec_zero-of-D_tq}) or (\ref{eq:spec_zero-of-D_aa}).
(Hence from Lemma \ref{lem:y-coincide},
 the polynomial representation is $Y$-semisimple.)
\end{enumerate}

Recall the definition of $c_{i,\lambda}$ in Proposition \ref{prop:int}.
We see that under the specialization of the parameters, 
 $c_{i,\lambda}^{\pm1}$ has no pole or zero
for any $0\leq i\leq n$ and $\lambda\in\mathbb{Z}^n$ such that
$\lambda\neq s_i\cdot\lambda$.
Therefore any non-symmetric Koornwinder polynomial $E_\lambda$ is cyclic
and the polynomial representation is irreducible for the specialized parameter.
\end{proof}

Hereafter, we will treat a specialization of parameters
 where the parameters are roots of one of the Laurent polynomials
(\ref{eq:spec_tq}), (\ref{eq:spec_aa}),
(\ref{eq:spec_ab}),(\ref{eq:spec_ac}), (\ref{eq:spec_ad}),
but the parameters are not roots of the other Laurent polynomials.
We call such a specialization {\it exclusive}.

\begin{defn}[exclusive specialization of parameters]\normalfont
We realize an exclusive specialization of parameters as follows.
Let $s\in \mathcal{A}$
 be one of the irreducible factors of (\ref{eq:spec_tq}), (\ref{eq:spec_aa}),
(\ref{eq:spec_ab}),(\ref{eq:spec_ac}), (\ref{eq:spec_ad}).
Then $\mathcal{A}/s\mathcal{A}$ is an integral domain and
 we denote the fractional field of $\mathcal{A}/s\mathcal{A}$ by $\mathbb{K}_s$.
We call $s$ a {\it specialization polynomial}.
\end{defn}

\begin{defn}[order of zeros or poles]\normalfont
Let $s$ be a specialization polynomial.
For any $a\in \mathcal{A}$, there exists $m\in\mathbb{Z}_{\geq0}$
such that $s^{-m}a\in \mathcal{A}\setminus s\mathcal{A}$.
We denote this integer $m$ by $\zeta_{s=0}(a)$, or simply by $\zeta(a)$.
For any $c=a/b\in \mathbb{K}$ $(a,b\in \mathcal{A})$,
define $\zeta(c)=\zeta(a)-\zeta(b)$.
For any $f=\sum_\lambda c_\lambda x^\lambda\in P_n$,
define $\zeta(f)=\min_\lambda\{\zeta(c_\lambda)\}$.
\end{defn}

For any $a\in \mathcal{A}$, let $\overline{a}$ be the quotient image
 in $\mathcal{A}/s\mathcal{A}$.
For any $c=a/b$ such that $\zeta(c)\geq0$,
define $c|_{s=0}=\overline{s^{-\zeta(b)}a}/\overline{s^{-\zeta(b)}b}$.
For any $f=\sum_\lambda c_\lambda x^\lambda\in P_n$ such that $\zeta(f)\geq0$,
define $f|_{s=0}=\sum_\lambda c_\lambda|_{s=0} x^\lambda$.

Let $\HH^s$ and $P_n^s$
 be the $(C^\vee_n,C_n)$-DAHA and its polynomial representation
over the field $\mathbb{K}_s$.
Hereafter,
unless stated otherwise,
we express the elements in $\mathbb{K}_s$, $\HH^s$, or $P_n^s$
by the specialization map $|_{s=0}$ (or mentioning ``at $s=0$").
Without the specialization map $|_{s=0}$ (or without mentioning ``at $s=0$"),
any scalars, operators, or Laurent polynomials
 are elements in $\mathbb{K}$, $\HH$, or $P_n$.

\subsection{Modified intertwiners and modified polynomials}
For a specialization polynomial $s$ and some $\lambda\in\mathbb{Z}^n$,
the multiplicity of $Y$-eigenvalue $y(\lambda)|_{s=0}$ in $P_n^s$
is possibly greater than $1$,
or $E_\lambda$ possibly has poles at $s=0$.
We can cancel the poles
 by taking linear combinations of non-symmetric Koornwinder polynomials
whose $Y$-eigenvalues are equal,
but there are many choices for such combinations.

In this subsection, we introduce a linear combination
 $\bar{E}_\lambda$ of the form
\bean
\bar{E}_\lambda=E_\lambda+\sum_{\mu}
m_{\lambda\mu}\frac{\chi_0^*(E_{\lambda})}{\chi_0^*(E_{\mu})}E_\mu
 \qquad(m_{\lambda\mu}\in\mathbb{Q}).\label{eq:mod-poly-form}
\eean
The sum runs over $\mu$ such that $y(\lambda)|_{s=0}=y(\mu)|_{s=0}$.
By taking suitable $m_{\lambda\mu}$,
 $(\bar{E}_\lambda)|_{s=0}$ has no pole at $s=0$
 (it may be a generalized $Y$-eigenvector).
We call such a combination (\ref{eq:mod-poly-form})
 a {\it modified polynomial}.
We will give a basis of the polynomial representation $P_n^s$
 in terms of $(\bar{E}_\lambda)|_{s=0}$.

\medskip

A sketch of construction of the basis is as follows.

For a given $\bar{E}_\lambda$,
 we define a modification $\bar{\phi}_i$ of the intertwiner $\phi_i$.
This $\bar{\phi}_i$ sends $\bar{E}_\lambda$ to
 $c_{i,\lambda}\bar{E}_{s_i\cdot\lambda}$.
We call it a {\it modified intertwiner}.
The modified intertwiner gives a inductive and combinatorial way
 to compute the coefficients $m_{\lambda\mu}$
 in the linear combination (\ref{eq:mod-poly-form}).

For any $\lambda\in\mathbb{Z}^n$, we take $w=s_{i_\ell}\cdots s_{i_1}\in W$
 such that $\lambda=w\cdot(0,\ldots,0)$
 and we obtain $\bar{E}_\lambda$ from $E_{(0,\ldots,0)}=1$
 by acting modified intertwiners
 $\bar{\phi}_{i_1}$, \ldots, $\bar{\phi}_{i_\ell}$.
We will see that
 $\bar{E}_\lambda$ is monic,
 and the terms in $\bar{E}_\lambda$ are lower than $x^\lambda$
 with respect to the ordering $\succ$.

We note that the modified intertwiners $\bar{\phi}_i$
 do not satisfy the braid relations
 though the original intertwiners $\phi_i$ satisfy them.
Therefore, the coefficients $m_{\lambda\mu}$ in $\bar{E}_\lambda$
(\ref{eq:mod-poly-form}) depend on reduced expressions of $w$ above.

\medskip

We state two lemmas.
(They are corollaries of Section 2.)

\begin{lem}\label{lem:ratio-of-chi}
For any $0\leq i\leq n$ and any $\lambda\in\mathbb{Z}^n$, we have
\bea
\frac{\chi_0^*(E_{s_i\cdot\lambda})}{\chi_0^*(E_{\lambda})}c_{i,\lambda}
&=&\left.
\frac{N_i(Y^{\alpha_i^\vee})}{D_i(Y^{\alpha_i^\vee})}
\right|_\lambda \quad(1\leq i\leq n-1)\\
\frac{\chi_0^*(E_{s_n\cdot\lambda})}{\chi_0^*(E_{\lambda})}c_{n,\lambda}
&=&\left.
\frac{N_n(Y^{\alpha_i^\vee})}{D_n(Y^{\alpha_i^\vee})}
\right|_\lambda \\
\frac{\chi_0^*(E_{s_0\cdot\lambda})}{\chi_0^*(E_{\lambda})}c_{0,\lambda}
&=&\left.
\frac{N_0(Y^{\alpha_i^\vee})}{D_0(Y^{\alpha_i^\vee})}\right|_\lambda.
\eea
\end{lem}

\begin{lem}\label{lem:T'}
For any $0\leq i\leq n$ and any $\lambda\in\mathbb{Z}^n$, we have
\bea
T_i'\stackrel{\mathrm{def}}{=}T_i-t^{1/2}&=&
\phi_i-\frac{N_i(Y^{\alpha_i^\vee})}{D_i(Y^{\alpha_i^\vee})}
 \quad(1\leq i\leq n-1)\\
T_n'\stackrel{\mathrm{def}}{=}T_n-t_n^{1/2}&=&
\phi_n-\frac{N_n(Y^{\alpha_i^\vee})}{D_n(Y^{\alpha_i^\vee})} \\
T_0'\stackrel{\mathrm{def}}{=}U_n-u_n^{1/2}&=&
\phi_0-\frac{N_0(Y^{\alpha_i^\vee})}{D_0(Y^{\alpha_i^\vee})}.
\eea
\end{lem}

\bigskip

Before giving an explicit definition of the modified intertwiners,
we will show an example.

\begin{ex}\normalfont
Suppose $n=4$.
Let $s$ be an irreducible factor in $t^2q-1$, and $\lambda:=(0,1,0,1)$.
We can easily check that $E_\lambda|_{s=0}$ is well-defined.
There are many paths to generate another $\mu\in\mathbb{Z}^n$
 from $\lambda=(0,1,0,1)$ by applying $s_i\in W$.
For instance,
\medskip

$\lambda=(0,1,0,1)$
 $\stackrel{s_3}{\rightarrow}$ $(0,1,1,0)$
 $\stackrel{s_1}{\rightarrow}$ $(1,0,1,0)$
 $\stackrel{s_2}{\rightarrow}$ $(1,1,0,0)=s_2s_1s_3\cdot\lambda$.

\medskip

In this example, we will construct polynomials
$\bar{E}_{s_3\cdot\lambda}$,
$\bar{E}_{s_1s_3\cdot\lambda}$,
$\bar{E}_{s_2s_1s_3\cdot\lambda}$
from the given polynomial $E_\lambda$.
We will see that they are of the form (\ref{eq:mod-poly-form}),
well-defined at $s=0$,
and generalized $Y$-eigenfunctions at $s=0$.

Since $(D_3(Y^{\alpha_3^\vee})|_\lambda)|_{s=0}=(t^2q-1)|_{s=0}=0$,
 the intertwiner $\phi_3$ is not well-defined at $s=0$.
Put $\bar{\phi}_3:=T_3'=T_3-t^{1/2}$. Then we see that
\bea
\bar{\phi}_3 E_\lambda&=&t^{1/2}E_{s_3\cdot\lambda}
-\frac{\chi_0^*(E_{s_3\cdot\lambda})}{\chi_0^*(E_{\lambda})}t^{1/2}E_\lambda \\
&=&t^{1/2}\bar{E}_{s_3\cdot\lambda} \\
\mbox{where } \bar{E}_{s_3\cdot\lambda}&:=&E_{s_3\cdot\lambda}
-\frac{\chi_0^*(E_{s_3\cdot\lambda})}{\chi_0^*(E_{\lambda})}E_\lambda.
\eea
Since $\bar{\phi}_3|_{s=0}$ and $E_\lambda|_{s=0}$ are well-defined,
$(\bar{E}_{s_3\cdot\lambda})|_{s=0}$ is also well-defined.
Moreover, from the commuting relations between $\bar{\phi}_3$ and $Y_i$,
we see that $(\bar{E}_{s_3\cdot\lambda})|_{s=0}$
 is a generalized $Y$-eigenvector.

Since $(D_1(Y^{\alpha_1^\vee})|_{s_3\cdot\lambda})|_{s=0}=(t^2q-1)|_{s=0}=0$,
 the intertwiner $\phi_1$ is not well-defined at $s=0$.
Put $\bar{\phi}_1:=T_1'=T_1-t^{1/2}$. Then we see that
\bean
\bar{\phi}_1\bar{E}_{s_3\cdot\lambda}
&=&T_1'\left(E_{s_3\cdot\lambda}
 -\frac{\chi_0^*(E_{s_3\cdot\lambda})}{\chi_0^*(E_{\lambda})}E_\lambda\right)
 \nonumber\\
&=&t^{1/2}\left(E_{s_1s_3\cdot\lambda}
 -\frac{\chi_0^*(E_{s_1s_3\cdot\lambda})}{\chi_0^*(E_{s_3\cdot\lambda})}
 E_{s_3\cdot\lambda}\right) 
 \nonumber\\
&&-\frac{\chi_0^*(E_{s_3\cdot\lambda})}{\chi_0^*(E_{\lambda})}
 t^{1/2}\left(E_{s_1\cdot\lambda}
 -\frac{\chi_0^*(E_{s_1\cdot\lambda})}{\chi_0^*(E_{\lambda})}
 E_{\lambda}\right).\label{eq:example_1010_1}
\eean
From the previous lemma, we have
\bean
\frac{\chi_0^*(E_{s_1\cdot\lambda})}{\chi_0^*(E_{\lambda})}t^{1/2}&=&
\frac{N_1(Y^{\alpha_1^\vee})}{D_1(Y^{\alpha_1^\vee})}|_\lambda,
 \label{eq:example_1010_2}\\
\frac{\chi_0^*(E_{s_1s_3\cdot\lambda})}{\chi_0^*(E_{s_3\cdot\lambda})}t^{1/2}&=&
\frac{N_1(Y^{\alpha_1^\vee})}{D_1(Y^{\alpha_1^\vee})}|_{s_3\cdot\lambda}
\label{eq:example_1010_3}.
\eean
Since $Y^{\alpha_1^\vee}|_{\lambda}=Y^{\alpha_1^\vee}|_{s_3\cdot\lambda}$,
we see that (\ref{eq:example_1010_2})$=$(\ref{eq:example_1010_3}).
Hence
\bea
(\ref{eq:example_1010_1})
&=&t^{1/2}\bar{E}_{s_1s_3\cdot\lambda}, \quad\mbox{where}\\
\bar{E}_{s_1s_3\cdot\lambda}&:=&E_{s_1s_3\cdot\lambda}
-\frac{\chi_0^*(E_{s_1s_3\cdot\lambda})}{\chi_0^*(E_{s_3\lambda})}E_{s_3\lambda}
-\frac{\chi_0^*(E_{s_1s_3\cdot\lambda})}{\chi_0^*(E_{s_1\lambda})}E_{s_1\lambda}
+\frac{\chi_0^*(E_{s_1s_3\cdot\lambda})}{\chi_0^*(E_{\lambda})}E_\lambda.
\eea
Since $\bar{\phi}_1|_{s=0}$ and $(\bar{E}_{s_3\cdot\lambda})|_{s=0}$
 are well-defined,
$(\bar{E}_{s_1s_3\cdot\lambda})|_{s=0}$ is also well-defined.
Moreover, from the commuting relations between $\bar{\phi}_1$ and $Y_i$,
we see that $(\bar{E}_{s_1s_3\cdot\lambda})|_{s=0}$
 is a generalized $Y$-eigenvector.

Let $\bar{\phi}_2\in\HH$ be as follows:
\bea
\bar{\phi}_2&=&\phi_2
\frac{D_2(Y^{\alpha_2^\vee})}{D_2(Y^{\alpha_2^\vee})|_{s_1s_3\cdot\lambda}} \\
&&\times\Large\fbox{$\left(
\frac{Y^{\alpha_2^\vee}-Y^{\alpha_2^\vee}|_\lambda}
{Y^{\alpha_2^\vee}|_{s_1s_3\cdot\lambda}-Y^{\alpha_2^\vee}|_\lambda}
-\frac{N_2(Y^{\alpha_2^\vee})|_{s_1s_3\cdot\lambda}}
{N_2(Y^{\alpha_2^\vee})|_\lambda}
\frac{Y^{\alpha_2^\vee}-Y^{\alpha_2^\vee}|_{s_1s_3\cdot\lambda}}
{Y^{\alpha_2^\vee}|_\lambda-Y^{\alpha_2^\vee}|_{s_1s_3\cdot\lambda}}
\right)$}.
\eea
Note that $(D_2(Y^{\alpha_i^\vee})|_{s_1s_3\cdot\lambda})|_{s=0}
=(tq-1)|_{s=0}\neq0$.
Although $N_2(Y^{\alpha_2^\vee})|_\lambda=t^{-1/2}(t^{-2}q^{-1}-1)$
and $Y^{\alpha_2^\vee}|_{s_1s_3\cdot\lambda}-Y^{\alpha_2^\vee}|_\lambda
=tq-t^{-3}q^{-1}$ have poles at $s=0$,
we see that these poles are cancelled out in $\bar{\phi}_2$.
Indeed,
\bea
&&\LARGE\fbox{$\cdots$}\\
&&=\frac{Y^{\alpha_2^\vee}-t^{-3}q^{-1}}
{tq-t^{-3}q^{-1}}
-\frac{t^2q-1}{t^{-2}q^{-1}-1}
\frac{Y^{\alpha_2^\vee}-tq}
{t^{-3}q^{-1}-tq} \\
&&\mbox{(by putting $Z:=t^2q$)},\\
&&=\frac{Y^{\alpha_2^\vee}-t^{-1}Z^{-1}}
{t^{-1}Z-t^{-1}Z^{-1}}
-\frac{Z-1}{Z^{-1}-1}
\frac{Y^{\alpha_2^\vee}-t^{-1}Z}
{t^{-1}Z^{-1}-t^{-1}Z}
\eea
and the factors $Z-1$ in the denominators above are cancelled out.
Thus we obtain that $\bar{\phi}_2|_{s=0}\in\HH^s$ is well-defined.
We have
\bea
&&\bar{\phi}_2\bar{E}_{s_1s_3\cdot\lambda} \\
&&=t^{1/2}E_{s_2s_1s_3\cdot\lambda}
-c_{2,\lambda}\frac{D_2(Y^{\alpha_2^\vee})|_\lambda}
{D_2(Y^{\alpha_2^\vee})|_{s_1s_3\cdot\lambda}}
\frac{N_2(Y^{\alpha_2^\vee})|_{s_1s_3\cdot\lambda}}
{N_2(Y^{\alpha_2^\vee})|_\lambda}
\frac{\chi_0^*(E_{s_1s_3\cdot\lambda})}{\chi_0^*(E_{\lambda})}
E_{s_2\cdot\lambda} \\
&&=t^{1/2}E_{s_2s_1s_3\cdot\lambda}
-c_{2,s_1s_3\cdot\lambda}
\frac{\chi_0^*(E_{s_2s_1s_3\cdot\lambda})}{\chi_0^*(E_{s_1s_3\cdot\lambda})}
\frac{\chi_0^*(E_{\lambda})}{\chi_0^*(E_{s_2\cdot\lambda})}
\frac{\chi_0^*(E_{s_1s_3\cdot\lambda})}{\chi_0^*(E_{\lambda})}
E_{s_2\cdot\lambda} \\
&&=t^{1/2}\bar{E}_{s_2s_1s_3\cdot\lambda} \\
&&\mbox{where }\bar{E}_{s_2s_1s_3\cdot\lambda}:=
E_{s_2s_1s_3\cdot\lambda}
-\frac{\chi_0^*(E_{s_2s_1s_3\cdot\lambda})}{\chi_0^*(E_{s_2\cdot\lambda})}
E_{s_2\cdot\lambda}.
\eea
Since $\bar{\phi}_2|_{s=0}$ and $\bar{E}_{s_1s_3\cdot\lambda}|_{s=0}$
 are well-defined,
we see that $(\bar{E}_{s_2s_1s_3\cdot\lambda})|_{s=0}$ is also well-defined.
Moreover, from the commuting relations between $\bar{\phi}_2$ and $Y_i$,
we see that $(\bar{E}_{s_2s_1s_3\cdot\lambda})|_{s=0}$
 is a generalized $Y$-eigenvector.
\hfill$\Box$
\end{ex}

\bigskip

In this example, we have given
the operators $\bar{\phi}_i\in\HH$
and linear combinations of non-symmetric Koornwinder polynomials
 of the form (\ref{eq:mod-poly-form}).
They are examples of {\it modified intertwiners} and {\it modified polynomials}
which will be defined below.
From now on,
we introduce a general setting.

Recall the definition of $D_i(Y^{\alpha_i^\vee})$ and $N_i(Y^{\alpha_i^\vee})$ 
(see (\ref{eq:def-of-D}) and (\ref{eq:def-of-N})).
We easily see that $D_i(Y^{\alpha_i^\vee})|_\lambda\neq0$ for any $\lambda$,
 and $\lambda=s_i\cdot\lambda
\Leftrightarrow c_{i,\lambda}=0
\Leftrightarrow N_i(Y^{\alpha_i^\vee})|_\lambda=0$.
Note that $D_i(Y^{\alpha_i^\vee})$ is the denominator of $\phi_i$.

Now we define an element $\bar{\phi}_i$ in $\HH$.

\begin{defn}[modified intertwiners]\label{defn:mod-int}\normalfont
Let $s$ be a specialization polynomial.
Fix $0\leq i\leq n$
 and $\lambda\in\mathbb{Z}^n$ such that $\lambda\neq s_i\cdot\lambda$.
Let $S_i=S_i(\lambda)$ be a finite set in $\mathbb{Z}^n$
 such that (i) $\lambda\in S_i$,
 (ii) $\mu\neq s_i\cdot\mu$ for any $\mu\in S_i$,
 (iii) $y(\lambda)|_{s=0}=y(\mu)|_{s=0}$ for any $\mu\in S_i$.
\quad For $\mu,\nu\in S_i$,
we write $\mu\sim\nu$ if $Y^{\alpha_i^\vee}|_\mu=Y^{\alpha_i^\vee}|_\nu$.
Let $\tilde{S_i}=S_i/\sim$,
 and $\tilde{\mu}$ be the equivalence class of $\mu\in S_i$.
For $\mu\in S_i$,
put
\bea
n_{\mu\lambda}&=&\left.
\frac{N_i(Y^{\alpha_i^\vee})|_\mu}{N_i(Y^{\alpha_i^\vee})|_\lambda}
\right|_{s=0} \\
d_{\lambda\mu}&=&\left.
\frac{D_i(Y^{\alpha_i^\vee})|_\lambda}{D_i(Y^{\alpha_i^\vee})|_\mu}
\right|_{s=0}.
\eea
(Well-definedness of $n_{\mu\lambda}$ and $d_{\lambda\mu}$
 is stated in Proposition \ref{prop:mod-int-welldef} below.)
Then, we denote the following element in $\HH$
by $\bar{\phi}_i(\lambda,\tilde{S_i})$:\\
If $(D_i(Y^{\alpha_i^\vee})|_\lambda)|_{s=0}\neq0$, then put
\bea
\bar{\phi}_i(\lambda,\tilde{S_i})&:=&
\phi_i
\frac{D_i(Y^{\alpha_i^\vee})}{D_i(Y^{\alpha_i^\vee})|_\lambda}
\sum_{\tilde{\mu}\in \tilde{S_i}}
n_{\mu\lambda}
\frac{N_i(Y^{\alpha_i^\vee})|_\lambda}{N_i(Y^{\alpha_i^\vee})|_\mu}
\prod_{\tilde{\nu}\in \tilde{S_i}\setminus\{\tilde{\mu}\}}
\frac{Y^{\alpha_i^\vee}-Y^{\alpha_i^\vee}|_\nu}
{Y^{\alpha_i^\vee}|_\mu-Y^{\alpha_i^\vee}|_\nu}
\eea
and if $(D_i(Y^{\alpha_i^\vee})|_\lambda)|_{s=0}=0$, then put
\bea
\bar{\phi}_i(\lambda,\tilde{S_i})&:=&
T_i'
\sum_{\tilde{\mu}\in \tilde{S_i}}
d_{\lambda\mu}
\frac{D_i(Y^{\alpha_i^\vee})|_\mu}{D_i(Y^{\alpha_i^\vee})|_\lambda}
\frac{N_i(Y^{\alpha_i^\vee})|_\lambda}{N_i(Y^{\alpha_i^\vee})|_\mu}
\prod_{\tilde{\nu}\in \tilde{S_i}\setminus\{\tilde{\mu}\}}
\frac{Y^{\alpha_i^\vee}-Y^{\alpha_i^\vee}|_\nu}
{Y^{\alpha_i^\vee}|_\mu-Y^{\alpha_i^\vee}|_\nu}.
\eea
We call $\bar{\phi}_i(\lambda,\tilde{S_i})$ a {\it modified intertwiner}.
\end{defn}

We sometimes write them as $\bar{\phi}_i$ for simplicity.
We will show well-definedness of $\bar{\phi}_i|_{s=0}$.

\begin{prop}[well-definedness]\label{prop:mod-int-welldef}
Let $s$ be a specialization polynomial.
Then $n_{\mu\lambda}=
\frac{N_i(Y^{\alpha_i^\vee})|_\mu}{N_i(Y^{\alpha_i^\vee})|_\lambda}|_{s=0}$
and
$d_{\lambda\mu}=
\frac{D_i(Y^{\alpha_i^\vee})|_\lambda}{D_i(Y^{\alpha_i^\vee})|_\mu}|_{s=0}$
are well-defined and belong to $\mathbb{Q}$.
(Hence $\bar{\phi}_i\in\HH$ is well-defined.)
The specialized modified intertwiner $(\bar{\phi}_i)|_{s=0}$
 is a well-defined element in $\HH^s$.
\end{prop}

In order to prove it, we use the following lemma.
\begin{lem}\label{lem:reduce}
Let $s$ be a specialization polynomial.
Then $s$ is an irreducible Laurent polynomial in $\mathcal{A}$ of the form
$s=z-\omega$ where $z\in\mathcal{A}$ is a monic Laurent monomial
 and $\omega$ is a primitive $\ell$-th root of unity for some $\ell$.
If $(z'-1)\in s\mathcal{A}$ for a monic Laurent monomial $z'\in\mathcal{A}$,
then $z'=z^{\ell m}$ for some $m\in\mathbb{Z}$.
\end{lem}

\begin{proof}[Proof of Proposition \ref{prop:mod-int-welldef}]

In this proof, for simplicity,
we write $N_i(Y^{\alpha_i^\vee})|_\lambda$ or
$D_i(Y^{\alpha_i^\vee})|_\lambda$
as $N_i(\lambda)$ or $D_i(\lambda)$ for any $\lambda\in\mathbb{Z}^n$.

From Lemma \ref{lem:reduce} above,
 $s=z-\omega$ is an irreducible Laurent polynomial
 for a monic Laurent monomial $z$
and a primitive $\ell$-th root of unity for some $\ell$.

First we show the well-definedness of $n_{\mu\lambda}$ for $i=0$.
(The proof for the other $i$ and $d_{\lambda\mu}$ are similar.)

Since $y(\lambda)|_{s=0}=y(\mu)|_{s=0}$ for any $\mu\in S_0$,
 it holds that $N_0(\lambda)|_{s=0}=N_0(\mu)|_{s=0}$.
Hence if $N_0(\lambda)|_{s=0}\neq0$, then $n_{\mu\lambda}=1$.
Suppose that $N_0(\lambda)|_{s=0}=0$.
We have
\bea
N_0(\lambda)&=&u_n^{1/2}(q^{1/2}y(\lambda)_1-q^{1/2}c')(q^{1/2}y(\lambda)_1-q^{1/2}d')\\
&=&qc'd'u_n^{1/2}(c^*y(\lambda)_1-1)(d^*y(\lambda)_1-1).
\eea
Hence a factor either $(c^*y(\lambda)_1-1)$ or $(d^*y(\lambda)_1-1)$ vanishes
and the other factor does not vanish at $s=0$.
Assume that $(c^*y(\lambda)_1-1)|_{s=0}=0$.
Since $c^*y(\lambda)_1$ is a monic Laurent monomial,
from Lemma \ref{lem:reduce}, we see that
$c^*y(\lambda)=z^{\ell m_1}$ for some $m_1\in\mathbb{Z}$.
Similarly, $c^*y(\mu)=z^{\ell m_2}$ for some $m_2\in\mathbb{Z}$.
Since $\lambda\neq s_0\cdot\lambda$,
we have $N_0(\lambda)\neq0$ and $m_1\neq0$.
Thus
\bean
n_{\mu\lambda}
\stackrel{\mathrm{def}}{=}\left.\frac{N_0(\mu)}{N_0(\lambda)}\right|_{s=0}
=\left.\frac{z^{\ell m_2}-1}{z^{\ell m_1}-1}\right|_{s=0}
=\frac{m_2}{m_1}\in\mathbb{Q}. \label{eq:n-is-given-by-ratio}
\eean

Next, we show the well-definedness for $\bar{\phi}_0$
for the case $D_0(\lambda)|_{s=0}=0$.
Proofs for the other operators are similar.

Take the representatives of $\tilde{S}_0$ by
 $\{\mu^{(1)}=\lambda,\mu^{(2)},\mu^{(3)},\ldots\}$.
Suppose $D_0(\lambda)|_{s=0}=0$.
Then $D_0(\mu^{(k)})|_{s=0}=0$ for any $k$.
Note that $D_0(\mu^{(k)})=qy(\mu^{(k)})_1^2-1$
and $qy(\mu^{(k)})_1^2$ is a monic Laurent monomial.
Moreover, the power of $qy(\mu^{(k)})_1^2$
with respect to the variables of $q^{1/2},t^{1/2},t_n^{1/2},t_0^{1/2}$
are even numbers.
Since $s=z-\omega$ is irreducible,
the powers of $z$ with respect to the variables
 $q^{1/2},t^{1/2},t_n^{1/2},t_0^{1/2},u_n^{1/2},u_0^{1/2}$
contain at least one odd number.
Hence from Lemma \ref{lem:reduce},
  $D_0(\mu^{(k)})=z^{-2\ell m_k}-1$ ($m_k\in\mathbb{Z}_{\neq0}$).
It means that $Y^{\alpha_0^\vee}|_{\mu^{(k)}}=z^{\ell m_k}$.
Thus $N_0(\mu^{(k)})|_{s=0}\neq0$ and $m_i\neq m_j$ (for any $i\neq j$).
Therefore
\bea
\bar{\phi}_0
&=&\sum_{1\leq k\leq |\tilde{S}_0|}
\frac{2m_1}{2m_k}\frac{z^{-2\ell m_k}-1}{z^{-2\ell m_1}-1}
\frac{N_0(z^{\ell m_1})}{N_0(z^{\ell m_k})}
\prod_{1\leq j\leq |\tilde{S}_0|,j\neq k}
\frac{Y^{\alpha_0^\vee}-z^{\ell m_j}}
{z^{\ell m_k}-z^{\ell m_j}}.
\eea
Note that $\zeta_{s=0}(N_0(z^{\ell m_k}))=0$.
Since $\bar{\phi}_0$ is a rational function with respect to $z^\ell$,
put $z^\ell=Z$.
We obtain well-definedness of $\bar{\phi}_0|_{s=0}$
if $\bar{\phi}_0$ does not have any poles at $Z=1$.
Reduce $\bar{\phi}_0$ to a common denominator and
take Taylor expansions  at $Z=1$ of the numerator and the common denominator.
Then we see that all factors $Z-1$ in the denominator are cancelled out.
\end{proof}

From Lemma \ref{lem:ratio-of-chi} and \ref{lem:T'},
 we easily obtain other expressions of the modified intertwiners:
If $(D_i(Y^{\alpha_i^\vee})|_\lambda)|_{s=0}\neq0$ then
\bean
\bar{\phi}_i(\lambda,\tilde{S}_i)&=&
\sum_{\tilde{\mu}\in \tilde{S}_i}
\left(n_{\mu\lambda}
\frac{\chi_0^*(E_{s_i\cdot\lambda})}{\chi_0^*(E_{\lambda})}
\frac{\chi_0^*(E_{\mu})}{\chi_0^*(E_{s_i\cdot\mu})}
\frac{c_{i,\lambda}}{c_{i,\mu}}
\phi_i\frac{D_i(Y^{\alpha_i^\vee})}{D_i(Y^{\alpha_i^\vee})|_\mu}
\right.\nonumber\\
&&\qquad\qquad\left.\prod_{\tilde{\nu}\in \tilde{S}_i\setminus\{\tilde{\mu}\}}
\frac{Y^{\alpha_i^\vee}-Y^{\alpha_i^\vee}|_\nu}
{Y^{\alpha_i^\vee}|_\mu-Y^{\alpha_i^\vee}|_\nu}\right),
\label{eq:mod-int-another-express_phi}
\eean
and if $(D_i(Y^{\alpha_i^\vee})|_\lambda)|_{s=0}=0$ then
\bean
\bar{\phi}_i(\lambda,\tilde{S}_i)&=&
\sum_{\tilde{\mu}\in \tilde{S}_i}
\left(d_{\lambda\mu}
\frac{\chi_0^*(E_{s_i\cdot\lambda})}{\chi_0^*(E_{\lambda})}
\frac{\chi_0^*(E_{\mu})}{\chi_0^*(E_{s_i\cdot\mu})}
\frac{c_{i,\lambda}}{c_{i,\mu}}
T_i' \right.\nonumber\\
&&\qquad\qquad
\left.\prod_{\tilde{\nu}\in \tilde{S}_i\setminus\{\tilde{\mu}\}}
\frac{Y^{\alpha_i^\vee}-Y^{\alpha_i^\vee}|_\nu}
{Y^{\alpha_i^\vee}|_\mu-Y^{\alpha_i^\vee}|_\nu}\right).
\label{eq:mod-int-another-express_T}
\eean

\begin{defn}[modified polynomials]\normalfont
For a finite set $S'$ satisfying
 $\lambda\not\in S'$ and $y(\lambda)|_{s=0}=y(\mu)|_{s=0}$ ($\forall\mu\in S'$),
denote the following sum
 by $\bar{E}_{\lambda,\{(m_{\lambda\mu},\mu);\mu\in S'\}}$:
\bea
\bar{E}_{\lambda,\{(m_{\lambda\mu},\mu);\mu\in S'\}}:=
E_\lambda+\sum_{\mu\in S'}
m_{\lambda\mu}\frac{\chi_0^*(E_{\lambda})}{\chi_0^*(E_{\mu})}E_\mu
 \qquad(m_{\lambda\mu}\in\mathbb{Q}).
\eea
We call it a {\it modified polynomial}.
We denote the data $\{(m_{\lambda\mu},\mu);\mu\in S'\}$
 of modification terms by $\mathrm{mt}(\lambda)$.
For simplicity, we sometimes write
$\bar{E}_{\lambda,\mathrm{mt}(\lambda)}$ as $\bar{E}_{\lambda}$.
\end{defn}

\begin{prop}[recurrence formula for modified polynomials]
\label{prop:recur-mod-poly}
For a finite set $S'$ satisfying
 $\lambda\not\in S'$ and $y(\lambda)|_{s=0}=y(\mu)|_{s=0}$ $(\forall\mu\in S')$,
suppose that $(\bar{E}_{\lambda,\{(m_{\lambda\mu},\mu);\mu\in S'\}})|_{s=0}$
is well-defined.
Fix $0\leq i\leq n$ and put $S_i:=S'\setminus\{\mu\in S';\mu=s_i\cdot\mu\}$
 and $\mathrm{mt}(\lambda):=\{(m_{\lambda\mu},\mu);\mu\in S'\}$.

(i)
If $\lambda\neq s_i\cdot\lambda$ and $D_i(\lambda)|_{s=0}\neq0$ then
\bea
\bar{\phi}_i(\lambda,\widetilde{S_i\cup \{\lambda\}})
\bar{E}_{\lambda,\mathrm{mt}(\lambda)}
&=&c_{i,\lambda}\bar{E}_{s_i\cdot\lambda,\mathrm{mt}(s_i\cdot\lambda)}\\
\mbox{where}\quad \mathrm{mt}(s_i\cdot\lambda)&:=&
\{(n_{\mu\lambda}m_{\lambda\mu},s_i\cdot\mu);\mu\in S_i\}
\eea
and the both sides are well-defined at $s=0$.

(ii)
If $\lambda\neq s_i\cdot\lambda$ and $D_i(\lambda)|_{s=0}=0$ then
\bea
\bar{\phi}_i(\lambda,\widetilde{S_i\cup \{\lambda\}})
\bar{E}_{\lambda,\mathrm{mt}(\lambda)}
&=&c_{i,\lambda}\bar{E}_{s_i\cdot\lambda,\mathrm{mt}(s_i\cdot\lambda)} \\
\mbox{where}\quad \mathrm{mt}(s_i\cdot\lambda)&:=&\{(-1,\lambda)\}\cup
\{(d_{\lambda\mu}m_{\lambda\mu},s_i\cdot\mu),
(-d_{\lambda\mu}m_{\lambda\mu},\mu);\mu\in S_i\}
\eea
and the both sides are well-defined at $s=0$.

(iii)
Suppose $\lambda=s_i\cdot\lambda$.
Then we see that $D_i(\lambda)|_{s=0}\neq0$.
For $\nu\in S'$, we have
\bea
\bar{\phi}_i(\nu,\tilde{S}_i)\bar{E}_{\lambda,\mathrm{mt}(\lambda)}
&=&c_{i,\nu}m_{\lambda\nu}\frac{\chi_0^*(E_{\lambda})}{\chi_0^*(E_{\nu})}
\bar{E}_{s_i\cdot\nu,\mathrm{mt}(s_i\cdot\nu)} \\
\mbox{where}\quad \mathrm{mt}(s_i\cdot\nu)&:=&
\{(n_{\mu\nu}m_{\lambda\nu}^{-1}m_{\lambda\mu},s_i\cdot\mu);\mu\in S_i\}
\eea
and the both sides are well-defined at $s=0$.
\end{prop}

\begin{proof}
We show (i).
The proofs for (ii) and (iii) are similar.
We see that
\bea
\phi_i\frac{D_i(Y^{\alpha_i^\vee})}{D_i(Y^{\alpha_i^\vee})|_\mu}
\left(\prod_{\tilde{\nu}\in \tilde{S}_i\setminus\{\tilde{\mu}\}}
\frac{Y^{\alpha_i^\vee}-Y^{\alpha_i^\vee}|_\nu}
{Y^{\alpha_i^\vee}|_\mu-Y^{\alpha_i^\vee}|_\nu}\right)E_\mu
&=&c_{i,\mu}E_{s_i\mu},\ \mbox{and} \\
\phi_i\frac{D_i(Y^{\alpha_i^\vee})}{D_i(Y^{\alpha_i^\vee})|_\mu}
\left(\prod_{\tilde{\nu}\in \tilde{S}_i\setminus\{\tilde{\mu}\}}
\frac{Y^{\alpha_i^\vee}-Y^{\alpha_i^\vee}|_\nu}
{Y^{\alpha_i^\vee}|_\mu-Y^{\alpha_i^\vee}|_\nu}\right)E_{\nu'}&=&0 \\
&&\mbox{(if $\tilde{\nu'}\in \tilde{S}_i\setminus\{\tilde{\mu}\}$ or
$\nu'=s_i\cdot\nu'$)}.
\eea
Hence from the expression (\ref{eq:mod-int-another-express_phi}),
\bea
&&\bar{\phi}_i(\lambda,\tilde{S}_i)
\left(E_\lambda+\sum_{\mu\in S'}
m_{\lambda\mu}\frac{\chi_0^*(E_{\lambda})}{\chi_0^*(E_{\mu})}E_\mu\right)\\
&&\qquad =c_{i,\lambda}E_{s_i\lambda}+\sum_{\mu\in S_i}
c_{i,\lambda}n_{\mu\lambda}m_{\lambda\mu}
\frac{\chi_0^*(E_{s_i\cdot\lambda})}{\chi_0^*(E_{s_i\cdot\mu})}E_{s_i\cdot\mu}
\eea
\end{proof}

\begin{rem}\label{rem:Chered}\normalfont
(i)
For given $\lambda\in\mathbb{Z}^n$,
 $\mathrm{mt}(\lambda)$, and $w\in W$,
from Proposition \ref{prop:recur-mod-poly},
we can compute the following data inductively:
$\mathrm{mt}(s_{i_1}\cdot\lambda)$,
$\mathrm{mt}(s_{i_2}s_{i_1}\cdot\lambda)$,
$\ldots$,
$\mathrm{mt}(s_{i_\ell}\cdots s_{i_1}\cdot\lambda)
=\mathrm{mt}(w\cdot\lambda)$
where $w=s_{i_\ell}\cdots s_{i_1}$ is a reduced expression.
However, in general, the result $\mathrm{mt}(w\cdot\lambda)$
depends on reduced expressions of $w$.
This fact corresponds to the fact that there are many choices
for taking a basis of a generalized $Y$-eigenspace
 whose dimension $d\geq2$
even except for a scalar multiple.
(Recall Proposition \ref{prop:int}.
the (non-modified) intertwiners $\phi_i$ satisfy the braid relations.
Hence $\phi_w:=\phi_{i_\ell}\cdots \phi_{i_1}$ is well-defined.
This fact corresponds to the fact that
if multiplicity of $Y$-eigenvalue is $1$,
then the choice of a basis of the $Y$-eigenspace
is unique up to a scalar multiple.)

(ii)
For the DAHA of reduced affine root systems,
Cherednik introduced chains of intertwiners 
and {\it non-semisimple Macdonald polynomials} in \cite{Ch}.
For $(C_n^\vee,C_n)$-DAHA,
by applying the intertwiners $\phi_i$ and the operators $T_i$,
one can also obtain
{\it non-semisimple non-symmetric Koornwinder polynomials}
(which are generalized $Y$-eigenfunctions).
However, coefficients of lower terms in the non-semisimple polynomials
 will become more complicated than the modified polynomials defined above.

On the contrary, in this paper,
we have defined rather complicated modified intertwiners.
However, from Proposition \ref{prop:recur-mod-poly},
we obtain the data $\mathrm{mt}(s_i\cdot\lambda)$
of lower terms from a given $\mathrm{mt}(\lambda)$
 by appending data $n_{\mu\lambda}$ or $d_{\lambda\mu}\in\mathbb{Q}$.
These rational numbers are nothing but ratios of powers of $Y$-eigenvalues
(recall (\ref{eq:n-is-given-by-ratio})).
Consequently, the coefficients of lower terms and its computation
become simpler.
\end{rem}

\begin{thm}[construction of a basis]\label{thm:construct-basis}
Let $\lambda\in\mathbb{Z}^n$.
Take the shortest element $w\in W$ such that
$\lambda=w\cdot(0,\ldots,0)$.
There exists a reduced expression $w=s_{j_\ell}\cdots s_{j_1}$
such that
$\langle \lambda^{(m)},\alpha_{j_{m}} \rangle>0$ $(1\leq m\leq \ell)$
where $\lambda^{(m)}:=s_{j_m}\cdots s_{j_1}\cdot (0,\ldots,0)$.
Determine $\bar{E}_{\lambda^{(m)}}$ recursively by
\bea
\bar{E}_{\lambda^{(m)}}
:=c_{j_m,\lambda^{(m-1)}}^{-1}\bar{\phi}_{j_m}\bar{E}_{\lambda^{(m-1)}}.
\eea
Then we have $\bar{E}_\lambda|_{s=0}$ is well-defined.
The set $\{\bar{E}_\lambda|_{s=0};\lambda\in\mathbb{Z}^n\}$
 constitutes a $\mathbb{K}_s$-basis of the polynomial representation $P_n^s$.

\end{thm}

\begin{proof}
By the well-definedness of $\bar{\phi}_{j_m}|_{s=0}$
 and $\bar{E}_{\lambda^{(m-1)}}|_{s=0}$,
we see that \\
 $(c_{j_m,\lambda^{(m-1)}}\bar{E}_{\lambda^{(m)}})|_{s=0}$
 is well-defined.
Since $\langle \lambda^{(m)},\alpha_{j_m}\rangle>0$,
 from Lemma \ref{lem:mod-poly-triangular} (see below),
$c_{j_m,\lambda^{(m-1)}}$ is a monic Laurent monomial.
Therefore \\
 $c_{j_m,\lambda^{(m-1)}}|_{s=0}\neq0$ and
 $\bar{E}_{\lambda^{(m)}}|_{s=0}$ is well-defined.

From Lemma \ref{lem:mod-poly-triangular},
the highest term in $\bar{E}_\lambda$ is $x^\lambda$.
Therefore $\{\bar{E}_\lambda|_{s=0};\lambda\in\mathbb{Z}^n\}$ are
linearly independent.

\end{proof}

\begin{lem}[triangularity]\label{lem:mod-poly-triangular}
Assume that $m_{\lambda\mu}\neq0$
for any $(m_{\lambda\mu},\mu)\in \mathrm{mt}(\lambda)$.
We call $\mathrm{mt}(\lambda)$ is {\it triangular} if
$\mu\prec\lambda$ for any $(m_{\lambda\mu},\mu)\in \mathrm{mt}(\lambda)$.
If $\langle s_i\cdot\lambda,\alpha_i\rangle>0$,
then we have

(i) $s_i\cdot\lambda\succ\lambda$,
 and $c_{i,\lambda}$ is a monic Laurent monomial,

(ii) if $\mathrm{mt}(\lambda)$ is triangular,
then $\mathrm{mt}(s_i\cdot\lambda)$
 given by Proposition \ref{prop:recur-mod-poly} is also triangular.
\end{lem}
\begin{proof}
(i) is clear.

We will show (ii).
Let $\mathrm{mt}(s_i\cdot\lambda)$ be the data given in
Proposition \ref{prop:recur-mod-poly}.
For any $(m_{s_i\cdot\lambda,\nu},\nu)\in\mathrm{mt}(s_i\cdot\lambda)$,
we see $\nu=\lambda$, $\mu$, or $s_i\cdot\mu$
 for some $(m_{\lambda\mu},\mu)\in \mathrm{mt}(\lambda)$
such that $s_i\mu\neq\mu$.
If $s_i\cdot\mu\prec\mu$, then
 $s_i\cdot\lambda\succ\lambda\succ\mu\succ s_i\cdot\mu$.
If $s_i\cdot\mu\succ\mu$, then
 $s_i\cdot\lambda\succ s_i\cdot\mu$ because $s_i\cdot\lambda\succ\lambda$ and
$\lambda\succ\mu$.
In both cases, we have
 $s_i\cdot\lambda\succ s_i\cdot\mu$ and $s_i\cdot\lambda\succ \mu$.
\end{proof}

\begin{rem}[well-definedness test]
\normalfont
Proposition \ref{prop:recur-mod-poly} and Theorem \ref{thm:construct-basis}
gives a sufficient condition for well-definedness of $E_\lambda|_{s=0}$
as follows:
We have constructed a $\mathbb{K}_s$-basis
 $\{\bar{E}_{\lambda,\mathrm{mt}(\lambda)}|_{s=0};\lambda\in\mathbb{Z}^n\}$
 of $P_n^s$.
The construction is given by inductive (and combinatorial) computation of data
$\mathrm{mt}(\lambda)$ from $\mathrm{mt}((0,\ldots,0))=\emptyset$.
If $\mathrm{mt}(\lambda)=\emptyset$ for some $\lambda\in\mathbb{Z}^n$, then
$\bar{E}_{\lambda,\mathrm{mt}(\lambda)}=E_\lambda$
and we obtain the well-definedness of $E_\lambda|_{s=0}$.
\end{rem}

\section{Polynomial representations for specialized parameters}

In this section, we examine irreducibility and $Y$-semisimplicity of
the polynomial representation under the exclusive specialization of parameters.

\subsection{Arrow relation}
\label{subsect:property-of-mod-poly}

We will investigate
 whether certain non-symmetric Koornwinder (or modified) polynomials
 are generated from other polynomials by actions of the (modified) intertwiners,
 or not.
In this subsection, we introduce a notation for this purpose.

\begin{defn}[arrow relation]\label{defn:arrow-rel}
For given elements $\lambda,\lambda'\in\mathbb{Z}^n$
 and data $\mathrm{mt}(\lambda),\mathrm{mt}(\lambda')$,
suppose that
two modified polynomials
 $\bar{E}_\lambda:=E_{\lambda,\mathrm{mt}(\lambda)}$ and
 $\bar{E}_{\lambda'}:=E_{\lambda',\mathrm{mt}(\lambda')}$
 are well-defined at $s=0$.
Then we write 
\bea
 (\lambda,\mathrm{mt}(\lambda))
\rightarrow(\lambda',\mathrm{mt}(\lambda'))
\eea
if there exist sequences
 $\lambda^{(0)}=\lambda, \lambda^{(1)},\ldots,\lambda^{(\ell)}=\lambda'$,
 and $i_1,\ldots,i_\ell$
 satisfying the following properties:
\bea
(\bar{\phi}_{i_m}\bar{E}_{\lambda^{(m-1)}})|_{s=0}
=(c_m\bar{E}_{\lambda^{(m)}})|_{s=0}\neq0
\eea
where $c_m\in\mathbb{K}$ and $\zeta(c_m)=0$ for any $1\leq m\leq \ell$.

We write
\bea
 (\lambda,\mathrm{mt}(\lambda))
\leftrightarrow(\lambda',\mathrm{mt}(\lambda'))
\eea
if two arrow relations hold:
\bea
(\lambda,\mathrm{mt}(\lambda))
\rightarrow(\lambda',\mathrm{mt}(\lambda')), \\
 (\lambda',\mathrm{mt}(\lambda'))
\rightarrow(\lambda,\mathrm{mt}(\lambda)).
\eea
\end{defn}

For example,
from Lemma \ref{lem:mod-poly-triangular},
if $\langle s_i\cdot\lambda,\alpha_i\rangle>0$
and $\mathrm{mt}(\lambda)$ is triangular,
then
$(\lambda,\mathrm{mt}(\lambda))
\rightarrow(s_i\cdot\lambda,\mathrm{mt}(s_i\cdot\lambda))$
where $\mathrm{mt}(s_i\cdot\lambda)$ is given in
Proposition \ref{prop:recur-mod-poly}.

As a corollary of the previous section,
we easily obtain the following arrow relation.
\begin{prop}\label{prop:int-on-mod-poly_good}
Suppose that $\zeta(N_i(\lambda))=0$, $\zeta(D_i(\lambda))=0$,
 and $\zeta(c_{i,\lambda})=0$.
Take a modified polynomial
 $\bar{E}_\lambda:=E_{\lambda,\{(m_{\lambda\mu},\mu)\}_\mu}$
which is well-defined at $s=0$.
Then for the data $\{(m_{\lambda\mu},s_i\cdot\mu)\}_\mu$
 which is given in Proposition \ref{prop:recur-mod-poly},
the modified polynomial
 $\bar{E}_{s_i\cdot\lambda}
:=E_{s_i\cdot\lambda,\{(m_{\lambda\mu},s_i\cdot\mu)\}_\mu}$
is well-defined at $s=0$, and
\bea
 (\lambda,\{m_{\lambda\mu},\mu\}_\mu)
\rightarrow(s_i\cdot\lambda,\{m_{\lambda\mu},s_i\cdot\mu\}_{\mu}).
\eea
\end{prop}


\subsection{The case $t^{k+1}q^{r-1}=1$ ($k+1\geq 2$)}\label{subsect:wheel}

In this subsection,
 we treat the exclusive specialization of the case (\ref{eq:spec_tq}).
Fix $n\geq k+1\geq 2$ and $r-1\geq 1$.
Let $s\in\mathcal{A}$ be an irreducible factor of
\bea
&&t^{(k+1)/m}q^{(r-1)/m}-\omega_m
\eea
where $m=GCD(k+1,r-1)$ and $\omega_m$ is a primitive root of unity.

The purpose of the subsection is
to realize a series of subrepresentations
$I_1^{(k,r)}\subset I_2^{(k,r)}\subset\cdots
 \subset I_{\lfloor\frac{n}{k+1}\rfloor}^{(k,r)}$
 in the polynomial representation $P_n^s$
as ideals defined by vanishing conditions.
For $I_1^{(k,r)}$,
we also show irreducibility and give a linear basis.

First we focus on $I_1^{(k,r)}$ and give a labeling set of its basis.

\begin{defn}[neighborhood, admissibility]\normalfont
\label{defn:neighborhood}
Fix integers $(a,b)$ ($n\geq a\geq2$, $b\geq1$) and $\lambda\in\mathbb{Z}^n$.
Put $w:=w_\lambda^+$.
The pair $(i,j)$ is called a {\it $(a,b)$-neighborhood} if
$(i,j)$ satisfies three conditions as follows:
\bea
&(i)& |\rho(\lambda)_i|-|\rho(\lambda)_j|=a-1, \\
&(ii)& |\lambda_i|-|\lambda_j|\leq b, \\
&(iii)& \mbox{if $|\lambda_i|-|\lambda_j|=b$, then}, \\
&&(iii-1)\quad
 \mbox{$(\sigma(\lambda)_i,\sigma(\lambda)_j)=(+,+)$ and $i>j$, or} \\
&&(iii-2)\quad
 \mbox{$(\sigma(\lambda)_i,\sigma(\lambda)_j)=(-,-)$ and $i<j$, or} \\
&&(iii-3)\quad
 \mbox{$(\sigma(\lambda)_i,\sigma(\lambda)_j)=(-,+)$}.
\eea
We denote the number of $(a,b)$-neighborhoods in $\lambda$
 by $\sharp^{(a,b)}(\lambda)$.
If $\sharp^{(a,b)}(\lambda)=0$, then $\lambda$
 is called {\it $(a,b)$-admissible}.
Mostly we consider the case $(a,b)=(k+1,r-1)$.
We sometimes omit $(a,b)$ (we call it {\it type}) if it is clear.
\end{defn}

\begin{rem}[another definition of neighborhood, admissibility]\normalfont
\label{rem:anotherdef}
For $\lambda\in\mathbb{Z}^n$,
take the shortest element $w_\lambda^0\in W_0$ such that
$w_\lambda^0\lambda\in\mathbb{Z}_{\geq0}^n$.
Define the map $\mathbb{Z}^n\rightarrow\mathbb{Z}_{\geq0}^n$:
$\lambda\mapsto\lambda^0=w_\lambda^0\lambda$.
Then $\sharp^{(a,b)}(\lambda)$ is equal to the {\it number of
neighborhoods of type $(a,b)$} in $\lambda^0$ in the sense of
 \cite{Ka2} (Definition 3.6).
Especially, $\lambda$ is admissible if and only if
 $\lambda^0$ has no neighborhoods of type $(a,b)$ in the sense of \cite{Ka2}.
We will reduce some combinatorial arguments to those in \cite{Ka2}.
\end{rem}

We will describe a basis of the irreducible subrepresentation
in terms of the non-symmetric Koornwinder polynomials $E_\lambda$
 for some $\lambda$.
The following proposition gives well-definedness of $E_\lambda$ at $s=0$.

\begin{prop}\label{prop:wheel_poly-welldef}
For any $\lambda$ such that $\sharp(\lambda)\leq1$,
there is no $\mu\in\mathbb{Z}^n$ such that $\mu\neq\lambda$ and
$y(\lambda)|_{s=0}=y(\mu)|_{s=0}$
 $($resp. $y^*(\lambda)|_{s=0}=y^*(\mu)|_{s=0}$$)$.
As a corollary, $E_\lambda$ $($resp. $E_\lambda^*$$)$ is well-defined
 at $s=0$.
\end{prop}
\begin{proof}
If $y(\lambda)=y(\mu)$ at $s=0$, then
there exist integers $m_1,\ldots,m_n\in\mathbb{Z}$ satisfying
(i) $\lambda_i=\mu_i+(r-1)m_i$,
(ii) $\rho(\lambda)_i=\rho(\mu)_i+(k+1)m_i$,
and (iii) $\sigma(\lambda)_i=\sigma(\mu)_i$.

From (iii), we see that $w^0_\lambda=w^0_\mu$.
Hence there exist $m_1',\ldots,m_n'\in\mathbb{Z}$
(i') $\lambda^0_i=\mu^0_i+(r-1)m_i'$,
(ii') $\rho(\lambda^0)_i=\rho(\mu^0)_i+(k+1)m_i'$,
and (iii') $\sigma(\lambda^0)_i=\sigma(\mu^0)_i=+$.
We also see that $\sharp(\lambda^0)=\sharp(\lambda)$.

For any element $\nu\in\mathbb{Z}^n_{\geq0}$,
the pair of integers $\rho(\nu)$ given in this paper
 and the pair of half-integers $\rho(\nu)$ given in \cite{Ka2} (\S2.3)
differ only by the total shift $\frac{n-1}{2}$.
In \cite{Ka2}, Lemma 4.13, by a combinatorial argument,
it was shown that
there is no $\mu^0\neq\lambda^0$ satisfying
 (i'), (ii'), and $\sharp(\lambda^0)\leq1$.

Therefore there is no $\mu\neq\lambda$ such that $y(\lambda)=y(\mu)$ at $s=0$.
The proof for the dual polynomial is similar.
\end{proof}

Now we give one of the main theorems in this subsection.
\begin{thm}\label{thm:wheel_irr}
The space
 $I_1^{(k,r)}{}'=\mathrm{span}_{\mathbb{K}_s}
 \{E_\lambda|_{s=0};\lambda$ is $admissible\}$
is an irreducible representation.
\end{thm}

\begin{proof}
By the definition, $I_1^{(k,r)}{}'$ is closed under the action of
 $Y_1,\ldots,Y_n$.
From Lemma \ref{lem:wheel_border} and Lemma \ref{lem:wheel_good-arrow} below,
$I_1^{(k,r)}{}'$ is closed under the action of $\phi_0,\ldots,\phi_n$.
Hence $I_1^{(k,r)}{}'$ is closed under the action of $\HH^s$.

Take any non-zero $v\in I_1^{(k,r)}{}'$.
Then there exists an admissible $\lambda$ and $h\in\HH^s$
 such that $E_\lambda=hv$.
Take another admissible $\mu$.
From Lemma \ref{lem:wheel_arrow_to0} and
 Lemma \ref{lem:wheel_arrow_0tolarge} below,
we have $\lambda\leftrightarrow\lambda^0$,
$\mu\leftrightarrow\mu^0$, and
$\lambda^0\leftrightarrow (nM,\ldots,2M,M)\leftrightarrow \mu^0$
for a large enough $M$.
Therefore any non-zero $v\in I_1^{(k,r)}{}'$ is cyclic.
\end{proof}

Four steps of the proof of Theorem \ref{thm:wheel_irr}
are given as follows.

\begin{lem}\label{lem:wheel_border}
If $\lambda$ is admissible and $s_i\cdot\lambda$ is not admissible
 for some $0\leq i\leq n$,
 then $(\phi_iE_\lambda)|_{s=0}=0$.
\end{lem}
\begin{proof}
From the assumption,
$\sharp(\lambda)=0$ and $\sharp(s_i\cdot\lambda)\geq1$.
Recall the defining condition (i), (ii), (iii-a), (iii-b), (iii-c)
 of neighborhood in Definition \ref{defn:neighborhood}.
If $(j_1,j_2)$ is a neighborhood in $s_i\cdot\lambda$ for $i=0$ or $n$,
then $(j_1,j_2)$ is also a neighborhood in $\lambda$.
Therefore we have $i\neq0$ and $i\neq n$.

For $1\leq i\leq n-1$,
if $(j_1,j_2)$ is a neighborhood in $s_i\cdot\lambda$,
then $(s_i(j_1),s_i(j_2))$ is also a neighborhood in $\lambda$ except
for the cases $(j_1,j_2)=(i,i+1)$ or $(i,i+1)$.
Since there is no neighborhood in $\lambda$,
the only neighborhood in $s_i\cdot\lambda$ is one of $(i,i+1)$ or $(i+1,i)$.

If $(i,i+1)$ is a neighborhood in $s_i\cdot\lambda$, then
$(\sigma(s_i\cdot\lambda)_i,\sigma(s_i\cdot\lambda)_{i+1})=(-,-)$.
If $(i+1,i)$ is a neighborhood in $s_i\cdot\lambda$, then
$(\sigma(s_i\cdot\lambda)_i,\sigma(s_i\cdot\lambda)_{i+1})=(+,+)$.
In both cases, $\langle \rho(s_i\cdot\lambda),\alpha_i\rangle=-k$ and
$\langle s_i\cdot\lambda,\alpha_i\rangle=-(r-1)$.
Hence $((t^{1/2}Y^{\alpha_i}-t^{-1/2})|_\lambda)|_{s=0}=0$,
and  $c_{i,\lambda}|_{s=0}=0$.
Thus $(\phi_iE_\lambda)|_{s=0}=(c_{i,\lambda}E_{s_i\cdot\lambda})|_{s=0}=0$.
\end{proof}

\begin{lem}\label{lem:wheel_good-arrow}
If $\lambda$ and $s_i\cdot\lambda$ are admissible
 for some $0\leq i\leq n$,
then $\lambda\leftrightarrow s_i\cdot\lambda$.
\end{lem}
\begin{proof}
It is easy to see that $\zeta(N_0(\lambda))=
\zeta(D_0(\lambda))=
\zeta(c_{0,\lambda})=0$
and $\zeta(D_n(\lambda))=0$.
We see that $\zeta(c_{n,\lambda})\neq0$ or
 $\zeta(N_n(\lambda))\neq0$ only if
\bean
&&\mbox{$\lambda_n=(r-1)m$ and and $\rho(\lambda)_n=(k+1)m$
 \quad for some $m\in\mathbb{Z}$}.\label{eq:wheel_good-arrow_proof1}
\eean
If (\ref{eq:wheel_good-arrow_proof1}) holds, then
there exists a $((k+1)m+1,(r-1)m)$-neighborhood in $\lambda$,
and consequently there exists a $(k+1,r-1)$-neighborhood in $\lambda$
(see \cite{Ka2}, Lemma 4.7.)
However it is inconsistent with the admissibility of $\lambda$.
Therefore from Proposition \ref{prop:int-on-mod-poly_good},
$\lambda\leftrightarrow s_i\cdot\lambda$ for $i=0$ and $n$.

For $1\leq i \leq n-1$, we see that $\zeta(N_i(\lambda))\neq0$ or
 $\zeta(D_i(\lambda))\neq0$ or $\zeta(c_{i,\lambda})\neq0$
only if
\bean
&&\mbox{$\langle\lambda,\alpha_i\rangle=(r-1)m$, and}\nonumber\\
&&\mbox{$\langle\rho(\lambda),\alpha_i\rangle=(k+1)m+d$, and}
\label{eq:wheel_good-arrow_proof2}\\
&&\mbox{$(\sigma(\lambda)_i,\sigma(\lambda)_{i+1})=(+,+)$ or $(-,-)$}\nonumber\\
&&\mbox{\qquad for some $m\in\mathbb{Z}$ and $d=-1,0,1$}.\nonumber
\eean
If (\ref{eq:wheel_good-arrow_proof2}) holds, then
there exists a $((k+1)m+d+1,(r-1)m)$-neighborhood in $\lambda$, and consequently
there exists a $(k+1,r-1)$-neighborhood in $\lambda$
(see \cite{Ka2}, Lemma 4.7.)
However it is inconsistent with the admissibility of $\lambda$.
Therefore from Proposition \ref{prop:int-on-mod-poly_good},
$\lambda\leftrightarrow s_i\cdot\lambda$ for $1\leq i \leq n-1$.
\end{proof}

\begin{lem}\label{lem:wheel_arrow_to0}
For any admissible $\lambda$, we have $\lambda\leftrightarrow \lambda^0$,
 where $\lambda^0$ is defined in Remark \ref{rem:anotherdef}.
\end{lem}
\begin{proof}
For a reduced expression of $w^0_\lambda=s_{i_l}\cdots s_{i_1}$,
put $\lambda^{(j)}=s_{i_j}\cdots s_{i_1}\cdot\lambda$ ($1\leq j \leq l$).
We can easily check that $\lambda^{(j)}$ is admissible
 for any $1\leq j \leq l$.
Hence from Lemma \ref{lem:wheel_good-arrow},
we have $\lambda^{(j-1)}\leftrightarrow\lambda^{(j)}$
 for any $1\leq j \leq l$.
\end{proof}

\begin{lem}\label{lem:wheel_arrow_0tolarge}
For an admissible $\lambda\in\mathbb{Z}^n_{\geq0}$ and a large enough $M$,
define $\lambda^{j,M}\in\mathbb{Z}^n_{\geq0}$ by
$\lambda^{j,M}_i=\lambda_i$ if $\rho(\lambda)_i<n-j$ and
$\lambda^{j,M}_i=(\rho(\lambda)_i+1)M$ if $\rho(\lambda)_i\geq n-j$.
Then we have $\lambda^{j,M}\leftrightarrow\lambda^{j+1,M}$.
\end{lem}
\begin{proof}
Take shortest element $w\in W$ such that $\lambda^{j+1,M}=w\lambda^{j,M}$
and take a reduced expression $w=s_{i_l}\cdots s_{i_1}$.
Put $\lambda^{(m)}=s_{i_l}\cdots s_{i_1}\lambda^{j,M}$.
Then we see that $\lambda^{(m)}$ is admissible for any $0\leq m \leq l$.
Hence from Lemma \ref{lem:wheel_good-arrow},
 $\lambda^{j,M}\leftrightarrow\lambda^{j+1,M}$.
\end{proof}

\bigskip

We will give a series of subrepresentations
$I_1^{(k,r)}\subset I_2^{(k,r)}\subset\cdots
 \subset I_{\lfloor\frac{n}{k+1}\rfloor}^{(k,r)}$
in terms of vanishing conditions for Laurent polynomials.
We will show that the irreducible representation $I_1^{(k,r)}{}'$
 coincides with $I_1^{(k,r)}$.

\begin{defn}[$m$-wheel condition]\normalfont
\label{defn:wheel_wheel-condition}
Let $\{i_1,\ldots,i_{k+1}\}$ be distinct indexes in $\{1,\ldots,n\}$
and $\{\sigma_1,\ldots,\sigma_{k+1}\}$ be a set of signs $+$ or $-$.
For $m\in\mathbb{Z}/(k+1)\mathbb{Z}$,
we denote by $\sigma_m i_m\mapsto\sigma_{m+1}i_{m+1}$
 the constraint $z_{i_m}^{\sigma_m}tq^{p_m}|_{s=0}=z_{i_{m+1}}^{\sigma_{m+1}}$
for some $p_m\in\mathbb{Z}$ satisfying (i) and (ii):
\bea
(i)\ p_m\geq 0 &&\\
(ii)\ p_m=0&\Rightarrow&
 \mbox{$(\sigma_m,\sigma_{m+1})=(+,+)$ and $i_m<i_{m+1}$, or}\\
&&\mbox{$(\sigma_m,\sigma_{m+1})=(+,-)$, or} \\
&&\mbox{$(\sigma_m,\sigma_{m+1})=(-,-)$ and $i_m>i_{m+1}$.} 
\eea
We call a closed cycle of the arrows ``$\mapsto$"
 with length $k+1$ a {\it wheel}.
(Since $(t^{k+1}q^{r-1}-1)|_{s=0}$, we see that $p_1+\cdots+p_{k+1}=r-1$.)
The $m$-wheel condition for $f$ is given as follows:
$f(z)=0$ if $z_1,\ldots,z_n$ form any disjoint $m$ wheels.
\end{defn}

Some examples for Laurent polynomials satisfying the $1$-wheel condition
 is given in \S\ref{subsect:example-wheel-condition}.

\begin{lem}\label{lem:wheel_alternative-wheel-condition}
The $1$-wheel condition for $f\in P_n^s$
 is equivalent to the following vanishing condition:
 $\chi_\mu^*(f)=0$ for any $\mu$ such that $\sharp(\mu)=1$.
\end{lem}
\begin{proof}
Take a neighborhood $(i,i')$ in $\mu$.
For $1\leq \ell\leq k+1$,
take $i_\ell$ and $\sigma_\ell$
 satisfying $|\rho(\mu)_{i_\ell}|=|\rho(\mu)_{i}|-\ell+1$
 and $\sigma_\ell=\sigma(\mu)_{i_\ell}$.
Then $\sigma_\ell i_\ell\mapsto\sigma_{\ell+1}i_{\ell+1}$ and
the cycle of the arrows forms a wheel.

Conversely, since $f$ is a Laurent polynomial,
 it is possible to replace the vanishing condition for
 the constraints $\sigma_m i_m\mapsto\sigma_{m+1}i_{m+1}$
by the vanishing condition for finitely many points in $\mathbb{K}_s^n$.
We can realize such points as
 $(\chi_\mu^*(x_1)|_{s=0},\ldots,\chi_\mu^*(x_n)|_{s=0})$
for some $\mu$ satisfying $\sharp(\mu)=1$.
\end{proof}

\begin{prop}
The space $I_m^{(k,r)}$ of Laurent polynomials
satisfying the $m$-wheel condition is a subrepresentation.
\end{prop}
\begin{proof}
Since the space $I_m^{(k,r)}$ is defined by the vanishing condition,
invariance under the action of $X_1^{\pm1},\ldots,X_n^{\pm1}$ is clear.

We show invariance under the action of $T_0,\ldots,T_n$.
Let $A$ and $B$ are coefficients of $s_i$ and $1$ in $T_i$: $T_i=As_i+B$.
Take $f\in I_m^{(k,r)}$
and fix a wheel: a set of constraints
 $w=\{\sigma_m i_m\mapsto\sigma_{m+1}i_{m+1}\}_m$.
Then $Bf$ vanishes under $w$.

The constraint $w$
for $s_if$ is given by
$s_i(z_{i_m})^{\sigma_m}tq^{p_m}|_{s=0}=s_i(z_{i_{m+1}})^{\sigma_{m+1}}$.
Hence it is equivalent to
$w'=\{s_i(\sigma_m i_m)\mapsto s_i(\sigma_{m+1}i_{m+1})\}_m$
provided that the powers of $q$ in the new constraints $w'$
satisfies the condition (i) and (ii)
 in Definition \ref{defn:wheel_wheel-condition}.
If (i) or (ii) is invalid, then we see that $A$ vanishes under $w$.
Hence $As_if$ vanishes under $w$.
\end{proof}

\medskip

Now we give the second main statement in this subsection.

\begin{thm}\label{thm:I1}
$I_1^{(k,r)}=I_1^{(k,r)}{}'$.
\end{thm}

The relation $I_1^{(k,r)}\supset I_1^{(k,r)}{}'$
is shown in Proposition \ref{prop:wheel_lower} below.

\begin{prop}\label{prop:wheel_lower}
If $\lambda$ is admissible, that is $\sharp(\lambda)=0$, then
 $\zeta(\chi_0^*(E_\lambda))=\lfloor\frac{n}{k+1}\rfloor$
 and $E_\lambda$ satisfies the wheel condition.
If $\lambda$ satisfies $\sharp(\lambda)=1$, then
 $\zeta(\chi_0(E_\lambda^*))=\lfloor\frac{n}{k+1}\rfloor-1$.
\end{prop}
\begin{proof}
We see that
\bea
\zeta(\chi_0^*(E_\lambda))&=&\zeta(\chi_0(E_\lambda^*)) \\
&=&\left\lfloor\frac{n}{k+1}\right\rfloor \\
&&-\sum_m\mbox{(the number of $((k+1)m,(r-1)m)$-neighborhoods in $\lambda$)} \\
&&+\sum_m\mbox{(the number of $((k+1)m+1,(r-1)m)$-neighborhoods in $\lambda$)}.
\eea
If $\lambda$ is admissible, then the second and the third term is zero.
If $\sharp(\lambda)=1$, then the second term is $1$ and the third term is zero.

Suppose that $\lambda$ is admissible and $\sharp(\mu)=1$.
By the duality relation, we have
\bea
\chi_\mu^*(E_\lambda)=\frac{\chi_\lambda(E_\mu^*)}{\chi_0(E_\mu^*)}
\chi_0^*(E_\lambda).
\eea
From Proposition \ref{prop:wheel_poly-welldef} and the first half of this proof,
we have $\zeta(\chi_\mu^*(E_\lambda))\geq1$.
Hence from the argument in the proof of
 Lemma \ref{lem:wheel_alternative-wheel-condition},
 we see that $E_\lambda$ satisfy the wheel condition.
\end{proof}

The rest of the proof of Theorem \ref{thm:I1} is 
 similar to Section 5 in \cite{Ka2}.
Consider a degree-restricted subspace of $I_1^{(k,r)}$,
introduce a dual space of the subspace.
By giving an ordering on monomials,
 obtain an upper bound of the dimension of the dual space.
We will see that the upper bound coincides with
 the lower bound which is given in Proposition \ref{prop:wheel_lower}.
\hfill $\Box$

\subsection{Examples of polynomials satisfying the wheel condition}
\label{subsect:example-wheel-condition}

In this subsection, we give two examples
 of {\it factorized} polynomials satisfying
 the $1$-wheel condition
(recall Definition \ref{defn:wheel_wheel-condition})
in the case $t^{k+1}q=1$ and $t^{k+1}q^k=1$.

Here we omit the specializing map $|_{s=0}$ for simplicity.
(All parameters are regarded as elements in $\mathbb{K}_s$.)

\begin{prop}\label{prop:level1sol}
Suppose that the parameters satisfy $t^{k+1}q=1$ $(n\geq k+1\geq 2)$
and that $n=km$.
Then
\bea
E_{(m-1,\ldots,1,0)^k}=
\prod_{\ell=1}^{k}\prod_{m(\ell-1)< i<j\leq m\ell}
 (x_i-t^{-1}x_j)(1-\frac{t^\ell q}{x_ix_j}).
\eea
\end{prop}

\begin{proof}
\underline{Step 1}\quad
Let us show that the RHS satisfies the wheel condition.

Consider a cycle of $k+1$ arrows.
Since $r-1=1$, the number of arrows from $-$ to $+$ is at most one.
Hence the following three cases are all of possibilities:
(i) all the signs are $+$,
(ii) all the signs are $-$,
(iii) the cycle is divided into one $+$ part and one $-$ part.
In the case (i) and (ii), the factor $(x_i-t^{-1}x_j)$ in the RHS vanishes.
Hereafter we assume (iii).
In this case, note that the power of $q$ in the arrow from $-$ to $+$ is $1$
and the powers of $q$ in the other arrows are $0$.

The indexes $\{1,\ldots,n\}$ are divided into $k$ blocks:
 $\{1,\ldots,m\},\ldots,\{n-m+1,\ldots,n\}$.
We enumerate these blocks from the left hand side.
(e.g. $\{n-m+1,\ldots,n\}$ is $k$-th block.)
In this proof we mean by $+\rightarrow+$ the arrow $+i\mapsto+j$ with $i<j$,
and by $-\leftarrow-$ the arrow $-i\mapsto-j$ with $i>j$.
If three vertices of the wheel are in the same block, then
there exists a pair of $+\rightarrow+$ or $-\leftarrow-$
and the factor $(x_i-t^{-1}x_j)$ vanishes.
So we assume there are at most two vertices in each block.

Let us show the following claim by induction:
``There exists $\ell$-th block such that
 $\ell$ blocks from the first to $\ell$-th blocks contain
 $\ell+1$ vertices of the wheel."

Consider the first block. By assumption, it contains at most two vertices.
If the number of vertices is two, then the claim holds.
We assume that it contains at most one vertex.

Suppose that $\ell$ blocks from the first to $\ell$-th blocks
contain at most $\ell$ vertices.
Then consider $\ell+1$-th block.
Since this block contains at most two vertices,
$\ell+1$ blocks from the first to $\ell+1$-th blocks
contain: (a) at most $\ell+1$ vertices, or
(b) $\ell+2$ vertices.
(b) implies the claim.
If we are in the case (a), consider the next block.

Note that the number of blocks are finite,
and that there are totally $k$ blocks and $k+1$ vertices.
So the case (b) does occur at least once,
and the claim is proved.

Consider the $\ell$-th block in the claim.
If $\ell$-th block contains $(+,+)$ or $(-,-)$, then
the factor $(x_i-t^{-1}x_j)$ vanishes.
The remaining possibility is that $\ell$-th block contains $+$ and $-$.
Denote these two indexes by $i_+$ and $i_-$.
From (iii) , we see that
 the sequence of arrows from $i_-$ to $i_+$ consists of
 one arrow from $-$ to $+$ and the other arrows are of the form
$-\leftarrow-$ or $+\rightarrow+$.
Moreover, from the claim,
the arrows from $i_-$ to $i_+$ intertwines
$\ell+1$ indexes in $\ell$ blocks from the first to $\ell$-th block.
Hence the factor $(1-\frac{t^lq}{x_{i_+}x_{i_-}})$ vanishes.

\underline{Step 2}\quad
Let $\lambda=(m-1,\ldots,1,0,\quad m-1,\ldots,1,0,\ \ldots,\ m-1,\ldots,1,0)$
(repeated by $k$ times)
and
\bea
S=\{\mbox{ admissible elements }\}
\cap\{\mu\in\mathbb{Z}^n;\mu\preceq\lambda\}.
\eea
Then we can check that $\sharp S=1$.

\underline{Step 3}\quad
From Step 1, the RHS is written of the form $\sum_{\mu\in S}c_\mu E_\mu$.
From Step 2, the RHS is equal to $c_\lambda E_\lambda$.
Since the coefficient of $x^\lambda$
in each side is $1$, we have $c_\lambda=1$.
%
\end{proof}

\begin{ex}
Suppose that the parameters satisfy $t^{k+1}q^{k}=1$ and that $n=km$.
Then, the product
\bea
\prod_{1\leq i<j\leq n}(x_ix_j^{-1}-t^{-1})(x_j-t^{-1}x_i^{-1})
\eea
satisfies the wheel condition.
\end{ex}
\begin{proof}
There exists at least one arrow whose power is $0$.
If the arrow is of the form $+\rightarrow+$ or $-\leftarrow-$,
then the factor $(x_ix_j^{-1}-t^{-1})$ vanishes.
If the arrow is from $+$ to $-$,
the factor $(x_j-t^{-1}x_i^{-1})$ vanishes.
\end{proof}


\subsection{Preliminaries for \S\ref{subsect:tq} and \S\ref{subsect:aa}}
\label{subsect:pre-irr}

In \S\ref{subsect:tq} and \S\ref{subsect:aa} below,
 we will show that the polynomial representation for specialized parameters
 $P_n^s$ is irreducible.
The proof is organized by the following three steps:

\underline{Step Irr-1}\quad
Define ``large enough" elements in $\mathbb{Z}^n$.
For any non-zero Laurent polynomial $f\in P_n^s$,
show that $hf=E_\lambda|_{s=0}$
for an element ${}^\exists h\in\HH^s$ and
 a large enough element $\lambda\in\mathbb{Z}^n$.

\underline{Step Irr-2}\quad
Find a ``specific element" $\lambda\in\mathbb{Z}^n$
 such that $\mu\leftrightarrow \lambda$
for any large enough element $\mu\in\mathbb{Z}^n$.

\underline{Step Irr-3}\quad
For the specific element $\lambda\in\mathbb{Z}^n$ in Step 2,
show that $\lambda\rightarrow (0,\ldots,0)$.
Since $E_{(0,\ldots,0)}=1$ is a cyclic vector in $P_n^s$,
we obtain irreducibility of $P_n^s$.

In Step Irr-3, we will use the following arrow relations.

\begin{prop}\label{prop:int-on-mod-poly_1to1}
Suppose that $\langle s_i\cdot\lambda,\alpha_i\rangle>0$,
$y(\lambda)|_{s=0}=y(\mu)|_{s=0}$,
$E_{s_i\cdot\lambda}$ and $E_{\lambda,(-1,\mu)}$
 and $E_\mu$ are well-defined at $s=0$,
$\zeta(D_i(s_i\cdot\lambda))=0$, $\zeta(N_i(s_i\cdot\lambda))=0$,
$\zeta(\chi_0^*(E_\lambda))-\zeta(\chi_0^*(E_\mu))=-1$,
and $\zeta(c_{i,\lambda})=1$.
Then we have
\bea
s_i\cdot\lambda\rightarrow\mu,
\eea
namely, $\bar{\phi}_iE_{s_i\cdot\lambda}|_{s=0}={}^\exists c E_\mu|_{s=0}$
 where $\zeta(c)=0$.
\end{prop}
\begin{proof}
We have
\bea
\bar{\phi}_iE_{s_i\cdot\lambda}&=&c_{i,\lambda} E_\lambda \\
&=& c_{i,\lambda}\left(
E_{\lambda,(-1,\mu)}+\frac{\chi_0^*(E_\lambda)}{\chi_0^*(E_\mu)}E_\mu\right).
\eea
The first term vanishes at $s=0$
and the second term is of the form $cE_\mu$ where $\zeta(c)=0$.
\end{proof}

\begin{prop}\label{prop:int-on-mod-poly_2to1}
Suppose that $\langle \lambda,\alpha_i\rangle=0$,
$\langle \mu,\alpha_i\rangle>0$,
$y(\lambda)|_{s=0}=y(\mu)|_{s=0}$,
$E_{\lambda,(-1,\mu)}$ and $E_{s_i\cdot\mu}$ are well-defined at $s=0$,
$\zeta(D_i(\lambda))=0$,
$\zeta(\chi_0^*(E_\lambda))-\zeta(\chi_0^*(E_\mu))=-1$,
and $\zeta(c_{i,\mu})=1$.
Then we have
\bea
(\lambda,(-1,\mu))\rightarrow s_i\cdot\mu,
\eea
namely,
 $\phi_iE_{\lambda,(-1,\mu)}|_{s=0}={}^\exists c E_{s_i\cdot\mu}|_{s=0}$
 where $\zeta(c)=0$.
\end{prop}
\begin{proof}
We have
\bea
\phi_iE_{\lambda,(-1,\mu)}
&=&-\frac{\chi_0^*(E_\lambda)}{\chi_0^*(E_\mu)}c_{i,\mu}E_{s_i\cdot\mu}.
\eea
Denote the coefficient of $E_{s_i\cdot\mu}$ by $c$, then $\zeta(c)=0$.
\end{proof}


We will show well-definedness of
 certain non-symmetric Koornwinder polynomials and modified polynomials at $s=0$
by the following argument.

Let $f\in P_n$ be a given Laurent polynomial.
Then $f|_{s=0}$ is well-defined (namely,
 $\zeta_{s=0}(f)\geq0$) if there exists a subset $N\subset\mathbb{Z}^n$
satisfying the following two conditions:
for a large enough integer $M\gg \max_i(|\deg_{x_i} f|)$,
\begin{description}
\item[(Grid-i)] $\{(\chi_\nu^*(x_1)|_{s=0},\ldots,\chi_\nu^*(x_n)|_{s=0});\nu\in N\}$
 contains distinct $M^n=M\times M\times\cdots\times M$ points
 in $\mathbb{K}_s^n$,
\item[(Grid-ii)] $\zeta_{s=0}(\chi_\nu^*(f))\geq0$ for any $\nu\in N$.
\end{description}

From this point of view, 
we will show the well-definedness of $E_\lambda$ as follows:
\begin{prop}
\label{prop:well-def-of-poly}
For a given element $\lambda\in\mathbb{Z}^n$,
suppose that there exists a subset $N\subset\mathbb{Z}^n$ satisfying
(Grid-i) for a large enough $M\gg \max_i|\lambda_i|$ and
\begin{description}
\item[(Grid-iii)] $E_\nu^*$ has no pole at $s=0$ for any $\nu\in N$,
\item[(Grid-iv)] $\zeta(\chi_0(E_\nu^*))=\zeta(\chi_0^*(E_\lambda))$
 for any $\nu\in N$.
\end{description}
Then $E_{\lambda}$ has no pole at $s=0$.
\end{prop}
\begin{proof}
For any $\nu\in N$, from the duality relation,
\bea
&&\chi_\nu^*(E_{\lambda})
=\frac{\chi_0^*(E_\lambda)}{\chi_0(E_\nu^*)} \chi_\lambda(E_\nu^*).
\eea
Thus $\zeta(\chi_\nu^*(E_{\lambda}))\geq0$.
It implies that $E_{\lambda}$ has no pole.
\end{proof}

Similarly, we will show the well-definedness of
 a modified polynomial as follows:
\begin{prop}
\label{prop:well-def-of-mod-poly}
For given elements $\lambda,\mu\in\mathbb{Z}^n$
suppose that
 $y(\lambda)|_{s=0}=y(\mu)|_{s=0}$,
 and that there exists a subset $N\subset\mathbb{Z}^n$ satisfying
(Grid-i) for a large enough $M\gg \max_i\{|\lambda_i|,|\mu_i|\}$,
 (Grid-iii), and
\begin{description}
\item[(Grid-v)] $\zeta(\chi_0(E_\nu^*))=\zeta(\chi_0^*(E_\lambda))+1$
 for any $\nu\in N$.
\end{description}
Then $E_{\lambda,(-1,\mu)}
=E_\lambda-\frac{\chi_0^*(E_\lambda)}{\chi_0^*(E_\mu)}E_\mu$
 has no pole at $s=0$.
\end{prop}
\begin{proof}
Take any $\nu\in N$.
Then
\bea
&&\chi_\nu^*(E_{\lambda,(-1,\mu)}) \\
&&=\chi_\nu^*(E_\lambda)
 -\frac{\chi_0^*(E_\lambda)}{\chi_0^*(E_\mu)}\chi_\nu^*(E_\mu) \\
&&=\chi_0^*(E_\lambda)\left(\frac{\chi_\nu^*(E_\lambda)}{\chi_0^*(E_\lambda)}
 -\frac{\chi_\nu^*(E_\mu)}{\chi_0^*(E_\mu)}\right) \\
&&=\frac{\chi_0^*(E_\lambda)}{\chi_0(E_\nu^*)}
 \left(\chi_\lambda(E_\nu^*)-\chi_\mu(E_\nu^*)\right) \qquad\mbox{(duality)}.
\eea
Since $y(\lambda)|_{s=0}=y(\mu)|_{s=0}$,
 we have $\zeta(\chi_\lambda(E_\nu^*)-\chi_\mu(E_\nu^*))\geq1$.
Thus $\chi_\nu^*(E_{\lambda,(-1,\mu)})$ has no pole at $s=0$.
It implies that $E_{\lambda,(-1,\mu)}$ has no pole at $s=0$.
\end{proof}

\subsection{The case $tq^{r-1}=1$}\label{subsect:tq}

Fix $r-1\geq1$ and let $s\in\mathcal{A}$ be an irreducible factor of
\bea
tq^{r-1}-1.
\eea

The main theorem of this subsection is as follows:

\begin{thm}\label{thm:tq_irred}
The polynomial representation $P_n^s$ is irreducible.
\end{thm}

We follow the sketch of proof given in \S\ref{subsect:pre-irr}.

\subsubsection{Step Irr-1}

\begin{defn}[large enough elements]\normalfont
In this subsection, $\lambda\in\mathbb{Z}^n$ is called {\it large enough}
 if $\lambda$ satisfies
\bea
\mbox{$\lambda=\lambda^+$ and $\lambda_i-\lambda_{i+1}\geq 2(r-1)$}.
\eea
\end{defn}

\begin{lem}\label{lem:tq_well-def-of-grids}
Take any large enough $\lambda$.
Then there are no $\mu\neq\lambda$
satisfying $y(\mu)=y(\lambda)$ or $y^*(\mu)=y^*(\lambda)$ at $s=0$.
As a corollary, $E_\lambda$ and $E^*_\lambda$ have no pole at $s=0$.
\end{lem}
\begin{proof}
Take $\mu$ satisfying $y(\mu)=y(\lambda)$ or $y^*(\mu)=y^*(\lambda)$ at $s=0$.
Then there exist integers $m_1,\ldots,m_n$ such that
$\mu_i=\lambda_i+(r-1)m_i$, $\rho(\mu)_i=\rho(\lambda)_i+m_i$, and
 $\sigma(\mu)_i=\sigma(\lambda)_i$.
However there is no such $\mu$ except for $\mu=\lambda$.
\end{proof}

\begin{lem}\label{lem:tq_any-to-grids}
Let $f$ be any non-zero Laurent polynomial.
Then there exists large enough $\lambda$ such that $E_\lambda|_{s=0}\in\HH^s f$.
\end{lem}
\begin{proof}
Let $f$ be any Laurent polynomial.
Take $\nu=\sum_{i=1}^n n_i\varpi_i$ for some $n_i\gg0$.
Then 
$x^\nu f=c_\lambda x^\lambda+\sum_{\mu\not\succeq\lambda} c_\mu x^\mu$
where
 $\lambda$ is large enough and $c_\lambda\neq0$.
Since $E_\lambda|_{s=0}$ is well-defined,
$x^\nu f=c_\lambda E_\lambda|_{s=0}+\sum_{\mu\not\succeq\lambda} c_\mu' x^\mu$.
From Lemma \ref{lem:tq_well-def-of-grids}, the dimension of
 the generalized $Y$-eigenspace with respect to the eigenvalue
$y(\lambda)$ is one.
Hence there exists an element $H\in\HH^s$ such that $Hx^\nu f=E_\lambda|_{s=0}$.

\end{proof}

\subsubsection{Step Irr-2}\quad

Here, we call $(2(n-1)(r-1),\ldots,4(r-1),2(r-1),0)\in\mathbb{Z}^n$
the {\it specific element}.

Using Proposition \ref{prop:int-on-mod-poly_good},
we can easily check that
 $(2(n-1)(r-1),\ldots,4(r-1),2(r-1),0)\leftrightarrow\lambda$
for any large enough $\lambda$.

\medskip

\subsubsection{Step Irr-3}\quad

Put
\bea
\lambda^{m,l}&:=&
\Bigl(\ 2(n-1)(r-1),\ldots,2(m+2)(r-1),\\
&&\qquad (m(r-1))^{l},2(m+1)(r-1),(m(r-1))^{m+1-l}\ \Bigl)
\eea
for $0\leq m\leq n-2$ and $0\leq l \leq m+1$, and
\bea
\lambda^{n-1,0}&:=&
\Bigl(\ (n-1)(r-1),\ldots, (n-1)(r-1) \ \Bigl).
\eea
Note that the specific element defined above is $\lambda^{0,0}$.

\begin{prop}
(i) For any $\lambda=\lambda^{m,l}$ with $0\leq m\leq n-2$ and $0\leq l \leq m$,
 or $(m,l)=(n-1,0)$, the non-symmetric Koornwinder polynomial
 $E_\lambda$ is well-defined at $s=0$.

(ii) $\lambda^{m,0}\rightarrow\lambda^{m,m}$

(iii) Put $\lambda=\lambda^{m,m+1}$ and $\mu=\lambda^{m+1,0}$
 with $0\leq m\leq n-2$.
 Then the modified polynomial
 $\bar{E}_{\lambda,(-1,\mu)}\stackrel{\mathrm{def}}{=}
E_\lambda-\frac{\chi_0^*(E_\lambda)}{\chi_0^*(E_\mu)}E_\mu$
 is well-defined at $s=0$.

(iv) ${\lambda^{m,m}}\rightarrow {\lambda^{m+1,0}}$.
\end{prop}
\begin{proof}
(i) From the evaluation formula (Proposition \ref{prop:chi0_formula})
 and recurrence relations (Lemma \ref{lem:chi0_recur}),
we have $\zeta(\chi_0^*(E_\lambda))=0$.
On the other hand, for any large enough $\nu$,
 we have $\zeta(\chi_0(E^*_\nu))=0$.
Then by putting
\bea
N=\{\nu\in\mathbb{Z}^n;\mbox{$\nu$ is large enough}\},
\eea
Proposition \ref{prop:well-def-of-poly} implies that
 $E_{\lambda}$ has no pole at $s=0$.

(ii) Use Proposition \ref{prop:int-on-mod-poly_good}.

(iii)
It is easy to see that $y(\lambda)|_{s=0}=y(\mu)|_{s=0}$.
From the evaluation formula (Proposition \ref{prop:chi0_formula})
 and recurrence relations (Lemma \ref{lem:chi0_recur}),
 we have $\zeta(\chi_0^*(E_\lambda))=-1$ and $\zeta(\chi_0^*(E_\mu))=0$.

Thus from Proposition \ref{prop:well-def-of-mod-poly},
we obtain the well-definedness of $\bar{E}_{\lambda,(-1,\mu)}$ at $s=0$.

(iv)
Use Proposition \ref{prop:int-on-mod-poly_1to1}.
\end{proof}

\begin{cor}
We have $\lambda^{0,0}\rightarrow\lambda^{n-1,0}
=((n-1)(r-1),\ldots,(n-1)(r-1))$.
\end{cor}

We will finish Step Irr-3.
\begin{prop}
We have $((n-1)(r-1),\ldots,(n-1)(r-1))\rightarrow(0,\ldots,0)$.
\end{prop}

\begin{proof}
In this proof, for simplicity,
we omit the factor $(r-1)$ in each component.
For instance, we regard
$(\alpha_1,\ldots,\alpha_n)$ as $((r-1)\alpha_1,\ldots,(r-1)\alpha_n)$.

We denote a data $(\lambda,\{(m_{\lambda\mu},\mu)\}_\mu)$ by
 $\lambda+\sum_\mu m_{\lambda\mu}\mu$.

For $m=1,\ldots,n-1$, we will show that
\bean
&&(m^n) \label{eq:lem:tq-irr-proof-1}\\
&&\rightarrow (m^{n-1},0)-(m^{n-m-1},(m-1)^m,m) \label{eq:lem:tq-irr-proof-2}\\
&&\rightarrow ((m-1),m^{n-m-1},(m-1)^m) \label{eq:lem:tq-irr-proof-3}\\
&&\rightarrow ((m-1)^n). \label{eq:lem:tq-irr-proof-4}
\eean

The step from (\ref{eq:lem:tq-irr-proof-1}) to
(\ref{eq:lem:tq-irr-proof-2}) is given by
 iterating the following arrow relations:
\bea
&&(m^n) \\
&&\leftrightarrow(m-1,m^{n-1}) \\
&&\leftrightarrow(m^{n-3},m-1,m,m) \\
&&\rightarrow(m^{n-2},m-1,m) \\
&&\rightarrow(m^{n-1},m-1)-(m^{n-2},m-1,m),
\eea
and for $l=1,\ldots,m-1$
\bea
&&(m^{n-1},l)-(m^{n-m+l-1},(m-1)^{m-l},m) \\
&&\leftrightarrow (l-1,m^{n-1})-(m-1,m^{n-m+l-1},(m-1)^{m-l}) \\
&&\leftrightarrow
 (m^{n-m+l-3},l-1,m^{m-l+2})-(m^{n-m+l-3},m-1,m^{2},(m-1)^{m-l}) \\
&&\rightarrow (m^{n-m+l-2},l-1,m^{m-l+1})-(m^{n-m+l-2},m-1,m,(m-1)^{m-l}) \\
&&\rightarrow
 \Big((m^{n-m+l-1},l-1,m^{m-l})-(m^{n-m+l-2},l-1,m^{m-l+1})\Big) \\
&& \quad -(m-l+1)\times
 \Big((m^{n-m+l-2},m,m-1,(m-1)^{m-l}) \\
&& \qquad\qquad\qquad\qquad\qquad-(m^{n-m+l-2},m-1,m,(m-1)^{m-l})\Big) \\
&&\rightarrow
 (m^{n-m+l},l-1,m^{m-l-1})-(m^{n-m+l-2},(m-1)^2,m,(m-1)^{m-l-1}) \\
&&\leftrightarrow (m^{n-1},l-1)-(m^{n-m+l-2},(m-1)^{m-l+1},m).
\eea

The step from (\ref{eq:lem:tq-irr-proof-2}) to
(\ref{eq:lem:tq-irr-proof-3}) is given as follows:
\bea
&& (m^{n-1},0)-(m^{n-m-1},(m-1)^m,m) \\
&& \rightarrow (m^{n-m-1},(m-1)^m,-m)
 \qquad \mbox{(Proposition \ref{prop:int-on-mod-poly_2to1})}\\
&&\leftrightarrow ((m-1),m^{n-m-1},(m-1)^m).
\eea

The step from (\ref{eq:lem:tq-irr-proof-3}) to
(\ref{eq:lem:tq-irr-proof-4})
 is given by iterating the following arrow relations:
for $l=m,\ldots,n-2$, we have
\bea
&&((m-1),m^{n-l-1},(m-1)^l) \\
&&\leftrightarrow (m^{n-l-3},(m-1),m,m,(m-1)^l) \\
&&\rightarrow (m^{n-l-2},(m-1),m,(m-1)^l) \\
&&\rightarrow
 (m^{n-l-2},m,(m-1),(m-1)^l)-(m^{n-l-2},(m-1),m,(m-1)^l) \\
&&\rightarrow (m^{n-l-2},(m-1),(m-1),m,(m-1)^{l-1})
 \qquad\mbox{(Proposition \ref{prop:int-on-mod-poly_2to1})} \\
&&\leftrightarrow (m^{n-l-2},(m-1)^{l+1},m) \\
&&\leftrightarrow ((m-1),m^{n-l-2},(m-1)^{l+1}).
\eea

By repeating from (\ref{eq:lem:tq-irr-proof-1}) to (\ref{eq:lem:tq-irr-proof-4}),
 we obtain the desired statement.
\end{proof}

We have shown $\lambda^{0,0} \rightarrow (0,\ldots,0)$,
and finished the proof of Theorem \ref{thm:tq_irred}.


\subsection{The case $t^{k+1}q^{r-1}a^*{}^2=1$}\label{subsect:aa}
Fix $2n-2\geq k+1\geq0$ and $r-1\geq 1$.
Let $s\in\mathcal{A}$ be an irreducible factor of
\bea
t^{k+1}q^{r-1}a^*{}^2-1.
\eea

The main theorem in this subsection is as follows.

\begin{thm}\label{thm:aa_irred}
The polynomial representation $P_n^s$ is irreducible.
\end{thm}

Similarly to \S\ref{subsect:tq},
the proof is given by the three steps.

\subsubsection{Step Irr-1}\quad

\begin{defn}[large enough elements]\normalfont
In this subsection, $\lambda\in\mathbb{Z}^n$ is called {\it large enough}
 if $\lambda_i\geq r-1$ for any $1\leq i\leq n$.
\end{defn}

\begin{lem}\label{lem:aa_well-def-of-grids}
Take any large enough $\lambda$.
Then there are no $\mu\neq\lambda$
satisfying $y(\mu)=y(\lambda)$ or $y^*(\mu)=y^*(\lambda)$ at $s=0$.
As a corollary, $E_\lambda$ and $E^*_\lambda$ have no pole at $s=0$.
\end{lem}
\begin{proof}
Take $\mu$ satisfying $y(\mu)=y(\lambda)$ or $y^*(\mu)=y^*(\lambda)$ at $s=0$.
Then there exist integers $m_1,\ldots,m_n$ such that
$\mu_i=\lambda_i+(r-1)m_i$, $\rho(\mu)_i=\rho(\lambda)_i+m_i$, and
 $\sigma(\mu)_i=\sigma(\lambda)_i$.
However there is no such $\mu$ except for $\mu=\lambda$.
\end{proof}

Similarly to the argument in Lemma \ref{lem:tq_any-to-grids},
for any Laurent polynomial $f\in P_n^s$,
there exists an element $\lambda\in\mathbb{Z}^n$
such that
 $\lambda_i\geq r-1$ (for any $1\leq i \leq n$) and $E_\lambda|_{s=0}\in\HH^s f$.

\subsubsection{Step Irr-2}\quad

Here, we call $(r-1,\ldots,r-1)\in\mathbb{Z}^n$
the {\it specific element}.
Using Proposition \ref{prop:int-on-mod-poly_good},
we can easily check that
 $\lambda\leftrightarrow(r-1,\ldots,r-1)$
for any large enough $\lambda$.

\subsubsection{Step Irr-3}\quad

We will show $(r-1,\ldots,r-1)\rightarrow(0,\ldots,0)$.
For simplicity, for any integer $m$ and natural number $\ell$, 
we denote by $m^\ell$ the sequence $m,\ldots,m$ of $\ell$-times repetition.

\begin{lem}\label{lem:aa_arrows}
Suppose $n-1\geq k+1\geq 1$.
Put $R:=(r-1)^{n-(k+3)}$ and
\bea
\lambda^{(1)}=(R,r-1,(-(r-1))^{k+2}),&& \mu^{(1)}=(R,r-1,0^{k+2}), \\
\lambda^{(2)}=(R,(-(r-1))^{k+2},-(r-1)),&& \mu^{(2)}=(R,0^{k+2},-(r-1)).
\eea
Then for $i=1,2$,
 we have

 (I) $y(\lambda^{(i)})|_{s=0}=y(\mu^{(i)})|_{s=0}$,

(II) $\zeta(\chi_0^*(E_{\lambda^{(i)}}))=-1$, $\zeta(\chi_0^*(E_{\mu^{(i)}}))
=\zeta(\chi_0^*(E_{s_n\lambda^{(1)}}))=\zeta(\chi_0^*(E_{s_{n-1}\mu^{(2)}}))=0$,

(III)
 $E_{s_n\lambda^{(1)}}\rightarrow E_{\mu^{(1)}}$ and
 $\bar{E}_{\lambda^{(2)},(-1,\mu^{(2)})}\rightarrow E_{s_{n-1}\mu^{(2)}}$.
\end{lem}
\begin{proof}
(I) is clear.
(II) is from Proposition \ref{prop:chi0_formula} and Lemma \ref{lem:chi0_recur}.

(III) 
Put $N=\{\nu\in\mathbb{Z}^n;\mbox{$\nu$ is large enough}\}$.
Then from Proposition \ref{prop:well-def-of-poly} and
Proposition \ref{prop:well-def-of-mod-poly},
we see that four polynomials
$E_{\mu^{(i)}}$, $\bar{E}_{\lambda^{(i)},(-1,\mu^{(i)})}$,
$E_{s_n\lambda^{(1)}}$, and $E_{s_{n-1}\mu^{(2)}}$ have no pole
at $s=0$.
Thus from Proposition \ref{prop:int-on-mod-poly_1to1},
 we have $E_{s_n\lambda^{(1)}}\rightarrow E_{\mu^{(1)}}$
and from Proposition \ref{prop:int-on-mod-poly_2to1},
we have
 $E_{\lambda^{(2)},(-1,\mu^{(2)})}\rightarrow E_{s_{n-1}\mu^{(2)}}$.
\end{proof}

\begin{lem}\label{lem:aa_arrows2}
Let $\lambda^{(1)}$ and $\mu^{(1)}$ as above.
Then $\mu^{(1)}\rightarrow (\lambda^{(1)},\{(-1,\mu^{(1)})\})$.
\end{lem}
\begin{proof}
From Theorem \ref{thm:construct-basis}, we have
\bea
&& \mu^{(1)} \rightarrow(\lambda^{(1)},\{(n_\nu,\nu);\nu\in {}^\exists S\}).
\eea
On the other hand, if $y(\nu)|_{s=0}=y(\lambda^{(1)})|_{s=0}$
for some $\nu\in\mathbb{Z}^n$,
 then $\nu=\lambda^{(1)}$ or $\nu=\mu^{(1)}$.
Since
$\zeta\left(\frac{\chi_0^*(E_{\lambda^{(1)}})}{\chi_0^*(E_{\mu^{(1)}})}\right)
=-1$,
 the only possibility of the modification data
 $\{(n_\nu,\nu);\nu\in S\}$ is
\bea
\{(n_\nu,\nu);\nu\in S\}=(-1,\mu^{(1)}).
\eea
\end{proof}

\begin{lem}\label{lem:aa_arrows3}
Suppose that $2n-2\geq k+1>n-1$.
For $i=3,\ldots,6$,
let $\lambda^{(i)}$, $\lambda^{(i)}{}'$, 
$\mu^{(i)}$, $\mu^{(i)}{}'$ be as follows:
\bea
\mbox{(for odd $k$)},&& \\
\lambda^{(3)}&=&((r-1)^{n-\frac{k+1}{2}-1},-(r-1),0,0^{\frac{k+1}{2}-1}), \\
\lambda^{(3)}{}'&=&((r-1)^{n-\frac{k+1}{2}-1},0,-(r-1),0^{\frac{k+1}{2}-1}), \\
\mu^{(3)}&=&((r-1)^{n-\frac{k+1}{2}-1},0,0,0^{\frac{k+1}{2}-1}), \\
\mbox{(for even $k$)},&& \\
\lambda^{(4)}&=&((r-1)^{n-\frac{k}{2}-2},0,-(r-1),0,0^{\frac{k}{2}-1}), \\
\lambda^{(4)}{}'&=&((r-1)^{n-\frac{k}{2}-2},0,0,-(r-1),0^{\frac{k}{2}-1}), \\
\mu^{(4)}&=&((r-1)^{n-\frac{k}{2}-2},-(r-1),0,0,0^{\frac{k}{2}-1}), \\
\mbox{for $2k\geq2l\geq k+3$},&& \\
\lambda^{(5)}&=&((r-1)^{n-l-1},0^{2l-k-3},0,0,-(r-1),0,0^{k-l}), \\
\lambda^{(5)}{}'&=&((r-1)^{n-l-1},0^{2l-k-3},0,0,0,-(r-1),0^{k-l}), \\
\mu^{(5)}&=&((r-1)^{n-l-1},0^{2l-k-3},0,-(r-1),0,0,0^{k-l}), \\
\lambda^{(6)}&=&((r-1)^{n-l-1},0^{2l-k-3},0,0,-(r-1),0,0^{k-l}), \\
\mu^{(6)}&=&((r-1)^{n-l-1},0^{2l-k-3},0,-(r-1),0,0,0^{k-l}), \\
\mu^{(6)}{}'&=&((r-1)^{n-l-1},0^{2l-k-3},-(r-1),0,0,0,0^{k-l}).
\eea
Then we have
 (I) $y(\lambda^{(i)})=y(\mu^{(i)})$ at $s=0$,

(II) $\zeta(\chi_0^*(E_{\lambda^{(i)}}))=-1$, $\zeta(\chi_0^*(E_{\mu^{(i)}}))
=\zeta(\chi_0^*(E_{\lambda^{(i)}{}'}))
=\zeta(\chi_0^*(E_{\mu^{(6)}{}'}))=0$,

 (III)
 $E_{\lambda^{(i)}{}'}\rightarrow E_{\mu^{(i)}}$ $(i=3,\ldots,5)$ and
 $\bar{E}_{\lambda^{(6)},(-1,\mu^{(6)})}\rightarrow E_{\mu^{(6)}{}'}$.
\end{lem}
\begin{proof}
(I) is clear.
(II) is from Proposition \ref{prop:chi0_formula} and Lemma \ref{lem:chi0_recur}.

(III) 
Put
 $N=\{\nu\in\mathbb{Z}^n;\mbox{$\nu$ is large enough}\}$.
Then from Proposition \ref{prop:well-def-of-poly} and
Proposition \ref{prop:well-def-of-mod-poly},
we see that four polynomials
$E_{\mu^{(i)}}$, $\bar{E}_{\lambda^{(i)},(-1,\mu^{(i)})}$,
$E_{\lambda^{(i)}{}'}$, and $E_{\mu^{(6)}{}'}$ have no pole
at $s=0$.
Thus from Proposition \ref{prop:int-on-mod-poly_1to1},
we have $E_{\lambda^{(i)}{}'}\rightarrow E_{\mu^{(i)}}$ $(i=3,\ldots,5)$
and from Proposition \ref{prop:int-on-mod-poly_2to1},
we have $\bar{E}_{\lambda^{(6)},(-1,\mu^{(6)})}\rightarrow E_{\mu^{(6)}{}'}$.
\end{proof}

Combining these lemmas, we obtain Step Irr-3.

\begin{proof}[Proof of Step Irr-3]
\underline{The case $n-2\geqq k+1\geqq 0$.}\quad
We denote the repetition $(r-1)^{n-(k+3)}$ by $R$.
\bean
&& (R,(r-1)^{k+3}) \nonumber\\
&& \leftrightarrow(R,r-1,(-(r-1))^{k+1},r-1) \label{eq:aa_irred_proof1}\\
&& \stackrel{\phi_n}{\rightarrow}(R,r-1,0^{k+2}) \label{eq:aa_irred_proof2}\\
&& \rightarrow(R,r-1,(-(r-1))^{k+2})-(R,r-1,0^{k+2}) \label{eq:aa_irred_proof3}
\\
&& \leftrightarrow(R,(-(r-1))^{k+2},-(r-1))-(R,0^{k+2},-(r-1)) \label{eq:aa_irred_proof4}\\
&& \stackrel{\phi_{n-1}}{\rightarrow}(R,0^{k+1},-(r-1),0) \label{eq:aa_irred_proof5}\\
&& \leftrightarrow(R,0^{k+1},0,0) \nonumber\\
&& \leftrightarrow(0,\ldots,0).\nonumber
\eean
The step from (\ref{eq:aa_irred_proof1}) to (\ref{eq:aa_irred_proof2}) is
from Lemma \ref{lem:aa_arrows}-(1).
The step from (\ref{eq:aa_irred_proof2}) to (\ref{eq:aa_irred_proof3}) is
from Lemma \ref{lem:aa_arrows2}.
The step from (\ref{eq:aa_irred_proof4}) to (\ref{eq:aa_irred_proof5}) is
from Lemma \ref{lem:aa_arrows}-(2).

\underline{The case $k+1=n-1$.}\quad The procedure above stops in the third step.

\underline{The case $2n-2\geqq k+1>n-1$.}\quad
We have
\bean
&& ((r-1)^n) \nonumber\\
&& \leftrightarrow ((r-1)^{n-[\frac{k}{2}]-1},0^{[\frac{k}{2}]+1}) \label{eq:aa_thm_proof1}\\
&& \rightarrow ((r-1)^{n-[\frac{k}{2}]-2},0^{[\frac{k}{2}]+2}) \label{eq:aa_thm_proof2}
\eean
The step from (\ref{eq:aa_thm_proof1}) to (\ref{eq:aa_thm_proof2}) is given
as follows:
If $k+1$ is even,
\bean
&&(\ref{eq:aa_thm_proof1}) \nonumber\\
&& \leftrightarrow ((r-1)^{n-[\frac{k}{2}]-2},0^{\frac{k+1}{2}},-(r-1)) \nonumber\\
&& \leftrightarrow ((r-1)^{n-[\frac{k}{2}]-2},0,-(r-1),0^{\frac{k+1}{2}-1}) \label{eq:aa_thm_proof3}\\
&&\rightarrow(\ref{eq:aa_thm_proof2}) \nonumber
\eean
The step from (\ref{eq:aa_thm_proof3}) to (\ref{eq:aa_thm_proof2})
is from Lemma \ref{lem:aa_arrows3}-(3).
If $k+1$ is odd,
\bean
&&(\ref{eq:aa_thm_proof1}) \nonumber\\
&& \leftrightarrow ((r-1)^{n-[\frac{k}{2}]-2},0^{\frac{k}{2}+1},-(r-1)) \nonumber\\
&& \leftrightarrow ((r-1)^{n-[\frac{k}{2}]-2},0,0,-(r-1),0^{\frac{k}{2}-1}) \label{eq:aa_thm_proof_4}\\
&& \rightarrow ((r-1)^{n-[\frac{k}{2}]-2},-(r-1),0,0,0^{\frac{k}{2}-1}) \label{eq:aa_thm_proof_5}\\
&&\leftrightarrow(\ref{eq:aa_thm_proof2}). \nonumber
\eean
The step from (\ref{eq:aa_thm_proof_4}) to (\ref{eq:aa_thm_proof_5})
is from Lemma \ref{lem:aa_arrows3}-(4).

We have $(\ref{eq:aa_thm_proof2})\rightarrow(0^n)$
by iterating the following steps:
for $2k\geq2l\geq k+3$
\bean
&& ((r-1)^{n-l},0^{l}) \nonumber\\
&& \leftrightarrow ((r-1)^{n-l-1},0^{l},-(r-1)) \nonumber\\
&& \leftrightarrow ((r-1)^{n-l-1},0^{2l-k-3},0,0,0,-(r-1),0^{k-l}) \label{eq:aa_thm_proof_l1}\\
&& \rightarrow ((r-1)^{n-l-1},0^{2l-k-3},0,-(r-1),0,0,0^{k-l}) \label{eq:aa_thm_proof_l2}\\
&& \rightarrow ((r-1)^{n-l-1},0^{2l-k-3},0,0,-(r-1),0,0^{k-l}) \nonumber\\
&& \quad -((r-1)^{n-l-1},0^{2l-k-3},0,-(r-1),0,0,0^{k-l}) \label{eq:aa_thm_proof_l3}\\
&& \rightarrow ((r-1)^{n-l-1},0^{2l-k-3},-(r-1),0,0,0,0^{k-l}) \label{eq:aa_thm_proof_l4}\\
&& \leftrightarrow ((r-1)^{n-l-1},0^{l+1}). \nonumber
\eean
The step from (\ref{eq:aa_thm_proof_l1}) to (\ref{eq:aa_thm_proof_l2})
is from Lemma \ref{lem:aa_arrows3}-(5).
The step from (\ref{eq:aa_thm_proof_l2}) to (\ref{eq:aa_thm_proof_l3})
is given by applying the case (ii) of Proposition \ref{prop:recur-mod-poly}.
The step from (\ref{eq:aa_thm_proof_l3}) to (\ref{eq:aa_thm_proof_l4})
is from Lemma \ref{lem:aa_arrows3}-(6).
\end{proof}


\subsection{The case $t^{n-i}q^{r-1}a^*b^*{}^{\pm1}=1$}\label{subsect:ab}

Fix $n\geq i\geq1$, $r-1\geq1$ and a sign $\pm1$.
We consider the following specialization of parameters:
\bean
t^{n-i}q^{r-1}a^*b^*{}^{\pm1}-1=0. \label{eq:spec_ab2}
\eean
Let $s=s_\pm\in\mathcal{A}$ be an irreducible factor of the left hand side.

The purpose of this subsection is to give a unique composition series
of the polynomial representation $P_n^s$,
and to give a characterization of
the irreducible subspace for the case $s=s_+$.

First we check that each $Y$-eigenspace is one-dimensional.
\begin{prop}
For any $\lambda\in\mathbb{Z}^n$,
 there is no $\mu\neq\lambda$ such that $y(\lambda)=y(\mu)$.
As a corollary, $E_\lambda|_{s=0}$ are well-defined.
\end{prop}
\begin{proof}
This is a corollary of Lemma \ref{lem:y-coincide}.
\end{proof}

\begin{thm}
The space
\bea
V=\mathrm{span}_{\mathbb{K}_s}\{ E_\lambda|_{s=0} ;
 \mbox{$\lambda^+_i>r-1$, or $\lambda^+_i=r-1$ and $\sigma(\lambda)_i=+1$}\}
\eea
is the unique irreducible subrepresentation of $P_n^s$.
\end{thm}
\begin{proof}
We see that $\zeta(N_j(\lambda))=\zeta(D_j(\lambda))=\zeta(c_{j,\lambda})=0$
except for the case $j=n$, $\lambda_n=\pm(r-1)$, $\rho(\lambda)_n=\pm(n-i)$.
Hence from Proposition \ref{prop:int-on-mod-poly_good},
 we have $\lambda\rightarrow s_j\lambda$
 except for the case $j=n$, $\lambda_n=\pm(r-1)$, $\rho(\lambda)_n=\pm(n-i)$.

If $j=n$, $\lambda_n=-(r-1)$, $\rho(\lambda)_n=-(n-i)$, then
$\zeta(D_n(\lambda))=\zeta(c_{n,\lambda})=0$.
Hence $\phi_nE_\lambda|_{s=0}=c_{n,\lambda}E_{s_n\lambda}|_{s=0}\neq0$.

If $j=n$, $\lambda_n=r-1$, $\rho(\lambda)_n=n-i$, then
$\zeta(D_n(\lambda))=0$ and $\zeta(c_{n,\lambda})=1$.
Hence $\phi_nE_\lambda|_{s=0}=c_{n,\lambda}E_{s_n\lambda}|_{s=0}=0$.

Therefore the space $V$ and the quotient space $P_n^s/V$ are irreducible.
\end{proof}

If $s=s_+$, there is a characterization of the
irreducible subspace in terms of a vanishing condition.
\begin{lem}
Let
\bea
S=\{ \lambda\in\mathbb{Z}^n ;
 \mbox{$\lambda^+_i>r-1$, or $\lambda^+_i=r-1$ and $\sigma(\lambda)_i=+1$}\}.
\eea
For any $\lambda\in S$, $\zeta_{s_+=0}(\chi_0^*(E_\lambda))=1$.
For any $\mu\not\in S$, $\zeta_{s_+=0}(\chi_0(E_\mu^*))=0$.
\end{lem}
\begin{proof}
Check the evaluation formula.
(Proposition \ref{prop:chi0_formula} and Lemma \ref{lem:chi0_recur}).
\end{proof}
\begin{prop}[characterization of $V$]
We have
\bea
V=\{f\in P_n^{s_+};
 \chi_\mu^*(f)|_{s_+=0}=0 \quad\mbox{for any $\mu\not\in S$}\}.
\eea
\end{prop}
\begin{proof}
Use the duality relation.
The proof is similar to that in Proposition \ref{prop:wheel_lower}.
\end{proof}


\subsection{The case $t^{n-i}q^{r-1}a^*c^*{}^{\pm1}=1$ or
$t^{n-i}q^{r-1}a^*d^*{}^{\pm1}=1$}\label{subsect:ac,ad}

Fix $1\leq i\leq n$, $r-1\geq1$, and a sign $\pm1$.
We consider the following specialization of parameters:
\bean
t^{n-i}q^{r-1-\theta(\pm1)}a^*c^*{}^{\pm1}-1&=&0,
 \quad\mbox{or}\label{eq:spec_ac2} \\
t^{n-i}q^{r-1-\theta(\pm1)}a^*d^*{}^{\pm1}-1&=&0, \label{eq:spec_ad2}
\eean
where $\theta(+1)=1$ and $\theta(-1)=0$.
Let $s=s_\pm\in\mathcal{A}$ be an irreducible factor of the left hand side.

The purpose of this subsection is similar to \S\ref{subsect:ab}.
We give a unique composition series of the polynomial representation $P_n^s$,
and give a characterization of
the irreducible subspace for the case $s=s_+$.
Since proofs are also similar, we omit the proofs.

\begin{prop}
For any $\lambda\in\mathbb{Z}^n$,
 there is no $\mu\neq\lambda$ such that $y(\lambda)=y(\mu)$.
As a corollary, $E_\lambda|_{s=0}$ are well-defined.
\end{prop}

\begin{prop}
The space
\bea
V=\mathrm{span}_{\mathbb{K}_s}\{ E_\lambda|_{s=0} ; 
\mbox{$\lambda^+_i>r-1$, or $\lambda^+_i=r-1$ and $\sigma(\lambda)_i=-1 $}\}
\eea
is the unique irreducible subrepresentation of $P_n^s$.
\end{prop}

Hereafter, we assume that $s=s_+$.

\begin{lem}
Let
\bea
S=\{ \lambda\in\mathbb{Z}^n ;
 \mbox{$\lambda^+_i>r-1$, or $\lambda^+_i=r-1$ and $\sigma(\lambda)_i=-1 $}\}.
\eea
For any $\lambda\in S$, $\zeta_{s_+=0}(\chi_0^*(E_\lambda))=1$.
For any $\mu\not\in S$, $\zeta_{s_+=0}(\chi_0(E_\mu^*))=0$.
\end{lem}

\begin{prop}[characterization]
We have
\bea
V=\{f\in P_n^{s_+};
 \chi_\mu^*(f)|_{s_+=0}=0 \quad\mbox{for any $\mu\not\in S$}\}.
\eea
\end{prop}

\medskip

\begin{rem}\normalfont
In \cite{vDSt}, Equation (3.6),
van Diejen and Stokman treated the specialization of parameters
of the following form:
\bea
t_lt_mt^{n-1}q^N=1 \quad(0\leq l\neq m\leq3).
\eea
The correspondence of their parameters
 and the present parameters is given by $t_0=a,t_1=b,t_2=c,t_3=d$.
If $l=0$, then their specializations correspond to
our specializations
 (\ref{eq:spec_ab2}), (\ref{eq:spec_ac2}), (\ref{eq:spec_ad2})
 where $i=1$ and the signs are plus.
Proposition 3.7 in \cite{vDSt} shows that
certain {\it symmetric} Koornwinder polynomials for specialized parameters
vanish at certain finite grid points.
Therefore, the characterization of $V$ in this paper
can be considered as a non-symmetric version and a generalization to $i\geq2$.
\end{rem}


\subsection{The case $q^{r-1}=1$}\label{subsect:q}

In this subsection, we consider the case (\ref{eq:spec_tq}) with $k+1=0$.
That is, suppose that $r-1\geq1$ and
let $s\in\mathcal{A}$ be an irreducible factor of
\bea
q-\omega_{r-1}
\eea
where $\omega_{r-1}$ is a primitive $(r-1)$-th root of unity.

The purpose of this subsection is
to show that $P_n^s$ is $Y$-semisimple and
to give infinitely-many subrepresentations $V_\mu$ of $P_n^s$
which are labelled by partitions $\mu$ of length $\leq n$.
(We only consider partitions of length $\leq n$
and we omit ``of length $\leq n$" here, for simplicity.)
The inclusion relation $V_\mu\supseteq V_{\mu'}$ holds
 if and only if $\mu\leq \mu'$
 where $\leq$ is the dominance ordering of partitions.
Any successive quotient of them is isomorphic to each other.

\medskip

For any fixed $\lambda\in\mathbb{Z}^n$,
there are infinitely-many $\mu\in\mathbb{Z}^n$
such that $y(\lambda)|_{s=0}=y(\mu)|_{s=0}$
since $q|_{s=0}$ is a root of unity.
However, we have the following statement.
\begin{prop}
For any $\lambda\in\mathbb{Z}^n$,
 $E_\lambda|_{s=0}$ is well-defined.
\end{prop}
\begin{proof}
Recall the recursive construction of a basis of $P_n^s$ in
 Theorem \ref{thm:construct-basis} and Proposition \ref{prop:recur-mod-poly}.
Since $D_i(\lambda)|_{s=0}\neq0$
 for any $0\leq i\leq n$ and any $\lambda\in\mathbb{Z}^n$,
 we see that the modification term $\mathrm{mt}(\lambda)=\emptyset$
 and $\bar{E}_\lambda=E_\lambda$.
Well-definedness of the modified polynomial $\bar{E}_\lambda$ at $s=0$
 implies the desired statement.
\end{proof}

We will introduce a labelling set of subrepresentations.

\begin{defn}[$(r-1)$-quotient]\label{defn:q_quotient}
\normalfont
For any $\lambda\in\mathbb{Z}^n$,
take the indexes $1\leq i_1,\ldots,i_n\leq n$ by
\bea
(i_1,\ldots,i_n)=(|w_\lambda^+(1)|,\ldots,|w_\lambda^+(n)|).
\eea
Put $i_{n+1}=n+1$, $\lambda_{i_{n+1}}=0$, and $\sigma(\lambda_{i_{n+1}})=+1$.
Take non-negative integers $p_m$ ($1\leq m\leq n$) by
\bea
p_m&=&\left\lfloor
\frac{|\lambda_{i_m}|-|\lambda_{i_{m+1}}|-\beta_m}{r-1}\right\rfloor,\\
\mbox{where $\beta_m=1$}&\mbox{if}&
\mbox{$(\sigma(\lambda)_{i_m},\sigma(\lambda)_{i_{m+1}})=(+,+)$ and $i_m>i_{m+1}$, or}\\
&&\mbox{$(\sigma(\lambda)_{i_m},\sigma(\lambda)_{i_{m+1}})=(-,-)$ and $i_m>i_{m+1}$, or}\\
&&\mbox{$(\sigma(\lambda)_{i_m},\sigma(\lambda)_{i_{m+1}})=(-,+)$}, \\
\mbox{and $\beta_m=0$}&&\mbox{otherwise}.
\eea
Then we call the partition
\bea
\lambda^{\mathrm{quot}}&:=&\sum_{i=1}^n p_i\varpi_i \\
&=&(\sum_{i=1}^{n}p_i,\sum_{i=2}^{n}p_i,\ldots,\sum_{i=n}^{n}p_i)
\eea
the {\it $(r-1)$-quotient} of $\lambda$.
We define the dominant element $\lambda^{\mathrm{std}}\in\mathbb{Z}^n$ by
\bea
\lambda^{\mathrm{std}}
:=((r-1)\lambda^{\mathrm{quot}}_1,\ldots,(r-1)\lambda^{\mathrm{quot}}_n).
\eea
Note that
 $(\lambda^{\mathrm{std}})^{\mathrm{quot}}=\lambda^{\mathrm{quot}}$.
\end{defn}

For example, let $n=4$, $r-1=3$ and $\lambda=(-3,0,-9,13)$.
Then $(i_1,i_2,i_3,i_4)=(3,4,2,1)$,
$(p_1,p_2,p_3,p_4)=(1,2,0,0)$,
$3$-quotient of $\lambda$ is $\lambda^{\mathrm{quot}}=(3,2,0,0)$,
and $\lambda^{\mathrm{std}}=(9,6,0,0)$.

For any partition $\mu$, put
\bea
\mathbb{Z}^n_\mu&:=&\{ \mbox{ the $(r-1)$-quotient of $\lambda$ is $\mu$ } \},\\
\mathbb{Z}^n_{\geq\mu}&:=&
\bigsqcup_{\mu':\mathrm{partition},\mu'\geq\mu}\mathbb{Z}^n_{\mu'}.
\eea
Then $\mathbb{Z}^n$ is decomposed as follows:
\bea
\mathbb{Z}^n=\bigsqcup_{\mu:\mathrm{partition}} \mathbb{Z}^n_\mu.
\eea
Define vector spaces as follows:
\bea
V_{\mu}&:=&\mathrm{span}_{\mathbb{K}_s}
 \{ E_\lambda|_{s=0} ; \lambda\in \mathbb{Z}^n_{\mu}\},\\
V_{\geq\mu}&:=&\mathrm{span}_{\mathbb{K}_s}
 \{ E_\lambda|_{s=0} ; \lambda\in \mathbb{Z}^n_{\geq\mu}\}.
\eea

Now we give the main theorem of this subsection.

\begin{thm}\label{thm:q_main}
(i) For any partition $\mu$ and any $\lambda\in\mathbb{Z}^n_\mu$,
$\HH^s E_\lambda|_{s=0}$ coincides with $V_{\geq\mu}$.
That is, $V_{\geq\mu}$ is a subrepresentation of $P_n^s$.

(ii) For any partitions $\mu$ and $\nu$,
 we have $V_{\geq\mu}\subseteq V_{\geq\nu}$ if and only if $\mu\geq\nu$.

(iii) For any partition $\mu$,
\bea
V'_\mu:=V_{\geq\mu}\left/\left(\sum_{\mu'>\mu}V_{\geq\mu'}\right)\right.
\eea
 is finite dimensional and irreducible.
As ${\mathbb{K}_s}$-vector spaces, $V'_\mu\cong V_\mu$.
For any $\lambda\neq\lambda'\in\mathbb{Z}^n_\mu$,
 we have $y(\lambda)|_{s=0}\neq y(\lambda')|_{s=0}$.
That is, the dimension of each $Y$-eigenspace in $V'_\mu$ is $1$.

(iv) For any partition $\mu$ and $\nu$,
 we have $V'_\mu\cong V'_\nu$
 as $\HH^s$-modules.
\end{thm}

\begin{proof}
We prove the desired statements by using lemmas given after the proof.

(i) We have
\bea
\HH^s E_\lambda|_{s=0}
&=&\sum_{\ell\geq0,0\leq i_1,\cdots,i_\ell\leq n}
 \mathbb{K}_s(\phi_{i_1}\cdots\phi_{i_\ell}E_\lambda)|_{s=0} \\
&=&\sum_{\nu^{\mathrm{quot}}\geq \mu} \mathbb{K}_s E_\nu|_{s=0}
 \quad\mbox{(from Lemma \ref{lem:q_arrows})} \\
&=&V_{\geq\mu}.
\eea

(ii)
We see that
\bea
&&V_{\geq\mu}\subseteq V_{\geq\nu} \\
&&\Leftrightarrow \mathbb{Z}^n_{\geq\mu}\subseteq \mathbb{Z}^n_{\geq\nu} \\
&&\Leftrightarrow \mu\geq\nu.
\eea

(iii) It is clear that $V'_\mu\cong V_\mu$ as ${\mathbb{K}_s}$-vector spaces.
We see that $V'_\mu$ is finite dimensional
because $\mathbb{Z}^n_\mu$ is a finite set.
Suppose that $y(\lambda)|_{s=0}=y(\lambda')|_{s=0}$
 for some $\lambda,\lambda'\in\mathbb{Z}^n_\mu$.
Then
$\rho(\lambda)=\rho(\lambda')$, $\sigma(\lambda)=\sigma(\lambda')$,
and $\lambda_i\equiv\lambda'_i$ mod $(r-1)$ for any $1\leq i\leq n$.
Since $\lambda$ and $\lambda'$ are elements in $\mathbb{Z}^n_\mu$,
$\lambda^{\mathrm{quot}}$ should be equal to $\lambda'^{\mathrm{quot}}$.
Thus $\lambda_i=\lambda'_i$ for any $1\leq i\leq n$.
Therefore, the irreducibility of $V'_\mu$ follows
  from (i), (ii), and Lemma \ref{lem:q_connected}.

(iv)
We show that $V'_\mu\cong V'_\nu$
for any partition $\mu$ and $\nu=(0,\ldots,0)$.
For any $\lambda\in\mathbb{Z}^n_\mu$,
we define $\lambda'\in\mathbb{Z}^n$ as follows:
let $i_1,\ldots,i_n$ and $p_1,\ldots,p_n$ be that of
 Definition \ref{defn:q_quotient}.
Put
\bean
\lambda'_{i_m}:=\lambda_{i_m}-\sgn(\lambda_{i_m})(r-1)(\sum_{j=m}^n p_j)
 \qquad (1\leq m \leq n). \label{eq:q_mainthm_isom}
\eean
(For example,
if $\lambda=(-3,0,-9,13)$, then $\lambda'=(-3,0,-3,4)$.)
We see that $\lambda'$ is an element in $\mathbb{Z}^n_{(0,\ldots,0)}$,
and the map $\lambda\mapsto\lambda'$
 gives the isomorphism $\mathbb{Z}^n_\mu\cong\mathbb{Z}^n_{(0,\ldots,0)}$
as finite sets.
By the definition of $\lambda'$ (\ref{eq:q_mainthm_isom}),
we have $y(\lambda)|_{s=0}=y(\lambda')|_{s=0}$.
Therefore, the actions of $Y_1,\ldots,Y_n$ and
$\phi_0,\ldots,\phi_n$ on $E_\lambda$
coincides with those on $E_{\lambda'}$,
and the map $E_\lambda\mapsto E_{\lambda'}$
extends to an isomorphism of $\HH^s$-modules.
\end{proof}

The following two lemmas
 are tools for the proof of Theorem \ref{thm:q_main}.

\begin{lem}\label{lem:q_arrows}
Fix $0\leq i\leq n$ and
 $\lambda\in\mathbb{Z}^n_\mu$ such that $s_i\cdot\lambda\neq\lambda$.
Then we have:

(i)
$c_{i,s_i\cdot\lambda}|_{s=0}=0$
 if and only if
$(s_i\cdot\lambda)^{\mathrm{quot}}=\lambda^{\mathrm{quot}}+\varpi_{j}$
for some $j$.

(ii)
$c_{i,\lambda}|_{s=0}=0$
 if and only if
$(s_i\cdot\lambda)^{\mathrm{quot}}=\lambda^{\mathrm{quot}}-\varpi_{j}$
for some $j$.

(iii)
$c_{i,\lambda}|_{s=0}\neq0$ and $c_{i,s_i\cdot\lambda}|_{s=0}\neq0$
 if and only if
 $(s_i\cdot\lambda)^{\mathrm{quot}}=\lambda^{\mathrm{quot}}$.

(The case $c_{i,\lambda}|_{s=0}=0$ and $c_{i,s_i\cdot\lambda}|_{s=0}=0$
does not occur by the definition of $c_{i,\lambda}$
 (see Proposition \ref{prop:int}).)
\end{lem}

\begin{cor}
We have the following facts.
(Recall the arrow relation ``$\rightarrow$"
 defined in Definition \ref{defn:arrow-rel}.)
In the case (i),
 we have $\lambda\rightarrow s_i\cdot\lambda$
 and $(\phi_iE_{s_i\cdot\lambda})|_{s=0}=0$.
In the case (ii), we have $s_i\cdot\lambda\rightarrow \lambda$
  and $(\phi_iE_{\lambda})|_{s=0}=0$.
In the case (iii), we have $\lambda\leftrightarrow s_i\cdot\lambda$.
\end{cor}

For example, if $n=4$, $r-1=3$ and $\lambda=(-3,0,-9,13)$, then
\bea
s_0&:& (-3,0,-9,13)\leftrightarrow(2,0,-9,13), \\
s_1&:& (-3,0,-9,13)\rightarrow(0,-3,-9,13), \\
s_2&:& (-3,0,-9,13)\leftrightarrow(-3,-9,0,13), \\
s_3&:& (-3,0,-9,13)\leftrightarrow(-3,0,13,-9), \\
s_4&:& (-3,0,-9,13)\leftrightarrow(-3,0,-9,-13).
\eea
$3$-quotient of $(0,-3,-9,13)$ is $(4,3,1,0)$,
which is equal to $(3,2,0,0)+(1,1,1,0)$.
$3$-quotient of the other elements are equal to $(3,2,0,0)$.

\begin{proof}[Proof of Lemma \ref{lem:q_arrows}]
First we give a proof for (i).
(The proof for (ii) is given by switching $s_i\cdot\lambda$ and $\lambda$.)

Recall the definition of $c_{i,\lambda}$ (see Proposition \ref{prop:int}).
Suppose that $c_{i,\lambda}|_{s=0}\neq0$ and $c_{i,s_i\cdot\lambda}|_{s=0}=0$.
Then $\langle s_i\cdot\lambda, \alpha_i\rangle >0$ and $1\leq i\leq n$.
Moreover, the condition (I) or (II) should be satisfied:

(I)
$1\leq i \leq n-1$,
$(\sigma(s_i\cdot\lambda)_i,\sigma(s_i\cdot\lambda)_{i+1})=(+1,+1)$ or $(-1,-1)$,
$\rho(s_i\cdot\lambda)_i-\rho(s_i\cdot\lambda)_{i+1}=1$, and
$(s_i\cdot\lambda)_i-(s_i\cdot\lambda)_{i+1}=(r-1)m$ for some $m>0$.

(II)
$i=n$,
$\sigma(s_i\cdot\lambda)_i=+1$,
$\rho(s_i\cdot\lambda)_i=0$, and
$(s_i\cdot\lambda)_i=(r-1)m$ for some $m>0$.

For each case (I) or (II),
 by the definition of $\lambda^{\mathrm{quot}}$,
we have $(s_i\cdot\lambda)^{\mathrm{quot}}=\lambda^{\mathrm{quot}}+\varpi_{j}$
for some $j$.

Conversely,
 if $(s_i\cdot\lambda)^{\mathrm{quot}}=\lambda^{\mathrm{quot}}+\varpi_{j}$
for some $j$, then
(I) or (II) occur, and we have
that $c_{i,\lambda}|_{s=0}\neq0$ and $c_{i,s_i\cdot\lambda}|_{s=0}=0$.

We show (iii).
By the definition of $\lambda^{\mathrm{quot}}$,
if $(s_i\cdot\lambda)^{\mathrm{quot}}\neq\lambda^{\mathrm{quot}}\pm\varpi_{j}$,
then $(s_i\cdot\lambda)^{\mathrm{quot}}=\lambda^{\mathrm{quot}}$.
Since we have proved (i) and (ii),
we obtain that $c_{i,\lambda}|_{s=0}\neq0$ and $c_{i,s_i\cdot\lambda}|_{s=0}\neq0$.
\end{proof}

\begin{lem}\label{lem:q_connected}
For any $\lambda\in\mathbb{Z}^n_\mu$,
 we have $\lambda\leftrightarrow\lambda^{\mathrm{std}}$.
\end{lem}
\begin{proof}
Fix $\lambda\in\mathbb{Z}^n_\mu$.
From the previous lemma, if $\lambda\neq s_i\cdot\lambda$ then
\bea
s_i\cdot\lambda\leftrightarrow\lambda \Leftrightarrow
(s_i\cdot\lambda)^{\mathrm{quot}}=\lambda^{\mathrm{quot}}.
\eea
Hence we will show $\lambda\leftrightarrow\lambda^{\mathrm{std}}$
 by applying simple reflections on $\lambda$ with preserving
 their $(r-1)$-quotients.
Such a procedure is realized as follows:

(Step 1)
In this step, we will change $\lambda\in\mathbb{Z}^n$
 to an element in $\mathbb{Z}^n_{\geq0}$ as follows.
Divide the set of indexes $\{1,\ldots,n\}$ into two pieces
$\{i_1<\cdots< i_\ell\}$
 such that $\lambda_{i_m}<0$ $(1\leq \forall i\leq \ell)$,
and
$\{j_1<\cdots< j_{n-\ell}\}$ such that $\lambda_{j_m}\geq0$
 $(1\leq \forall i\leq n-\ell)$.
Then putting
\bea
\lambda':=
(-\lambda_{i_\ell}-1,\ldots,-\lambda_{i_1}-1,
\lambda_{j_1},\ldots,\lambda_{j_{n-\ell}}),
\eea
we have $\lambda'^{\mathrm{quot}}=\lambda^{\mathrm{quot}}$
 and $\lambda'\leftrightarrow\lambda$.
(Indeed, for a shortest element $s_{k_L}\cdots s_{k_1}\in W$ such that
$\lambda'=s_{k_L}\cdots s_{k_1}\cdot\lambda$,
 the $(r-1)$-quotients of
 $\lambda$, $s_{k_1}\cdot\lambda$, $s_{k_2}s_{k_1}\cdot\lambda$,
 $\ldots$, $s_{k_L}\cdots s_{k_1}\cdot\lambda$ are equal.)

Note that $\lambda'\in\mathbb{Z}^n_{\geq0}$.
Hereafter we assume that $\lambda\in\mathbb{Z}^n_{\geq0}$.

(Step 2)
In this step, we will reduce $\lambda\in\mathbb{Z}^n_{\geq0}$
 to a dominant element $\lambda''$ as follows.
Take the index
\bea
(i_1,\ldots,i_n)=(|w_\lambda^+(1)|,\ldots,|w_\lambda^+(n)|).
\eea
Suppose that $i_{\ell+1}=\ell+1,i_{\ell+2}=\ell+2,\ldots,i_{n}=n$.
Then for any $j$ satisfying $i_\ell+1\leq j\leq \ell$,
we see that $\lambda_{i_\ell}<\lambda_{j}$.
Hence by putting
\bea
\lambda':=(\lambda_{i_\ell+1}-1,\ldots,\lambda_{\ell}-1,
\lambda_1,\lambda_2,\ldots,\lambda_{i_\ell},
\lambda_{i_{\ell+1}},\ldots,\lambda_{i_n}),
\eea
we have $\lambda'^{\mathrm{quot}}=\lambda^{\mathrm{quot}}$
 and $\lambda'\leftrightarrow\lambda$.
(The reason is similar to that in Step 1.)
Moreover new indexes
\bea
(i_1',\ldots,i_n'):=(|w_{\lambda'}^+(1)|,\ldots,|w_{\lambda'}^+(n)|)
\eea
satisfy $i_{\ell}'=\ell,i_{\ell+1}'=\ell+1,\ldots,i_{n}'=n$.
Therefore by repeating this inductively,
we obtain a dominant element $\lambda''$
such that $\lambda''^{\mathrm{quot}}=\lambda^{\mathrm{quot}}$
 and $\lambda''\leftrightarrow\lambda$.

Hereafter we assume that $\lambda$ is dominant.

(Step 3)
In this step, we will reduce a dominant element $\lambda\in\mathbb{Z}^n$
 to $\lambda^{\mathrm{std}}$ as follows.
Suppose that $\lambda_i=\lambda^{\mathrm{std}}_i$ for any $\ell+1\leq i \leq n$.
By the definition of $\lambda^{\mathrm{std}}$,
we see $\lambda_\ell\geq\lambda^{\mathrm{std}}_\ell$.
If $\lambda_\ell>\lambda^{\mathrm{std}}_\ell$, then
by putting
\bea
\lambda':=(\lambda_1-1,\ldots,\lambda_\ell-1,\lambda_{\ell+1},\cdots,\lambda_n),
\eea
we have $\lambda'^{\mathrm{quot}}=\lambda^{\mathrm{quot}}$
 and $\lambda'\leftrightarrow\lambda$.
(The reason is similar to that in (i).)
Therefore by repeating this inductively,
we obtain the desired element $\lambda^{\mathrm{std}}$
and it satisfies
 $(\lambda^{\mathrm{std}})^{\mathrm{quot}}=\lambda^{\mathrm{quot}}$
 and $\lambda^{\mathrm{std}}\leftrightarrow\lambda$.
\end{proof}

For example, let $n=4$, $r-1=3$ and $\lambda=(-3,0,-9,13)$.
Then from Step 1, we obtain $\lambda\leftrightarrow(8,2,0,13)$.
From Step 2, we obtain $(8,2,0,13)\leftrightarrow(12,8,2,0)$.
From Step 3, we obtain
 $(12,8,2,0)\leftrightarrow(9,6,0,0)=\lambda^{\mathrm{std}}$.
$3$-quotient of these elements is $(3,2,0,0)$.



\medskip


\begin{thebibliography}{9999}

\bibitem[Ch]{Ch}
I. Cherednik,
Non-semisimple Macdonald polynomials,
arXiv:0709.1742.

\bibitem[ChBook]{ChBook}
I. Cherednik,
Double affine Hecke algebras,
London Mathematical Society Lecture Note Series.

\bibitem[vDSt]{vDSt}
J.F. van Diejen, J.V. Stokman,
Multivariable $q$-Racah polynomials,
Duke Math. J. 91 (1998), 89--136.

\bibitem[En]{En}
N. Enomoto,
Composition Factors of Polynomial Representation of DAHA and Crystallized
Decomposition Numbers,
math.RT/0604368.

\bibitem[FJMM]{FJMM}
B. Feigin, M. Jimbo, T. Miwa, E. Mukhin,
Symmetric polynomials vanishing on the shifted diagonals and Macdonald polynomials,
Int. Math. Res. Not. 2003, no. 18, 1015--1034.

\bibitem[Ka1]{Ka1}
M. Kasatani,
Zeros of symmetric Laurent polynomials of type $(BC)\sb n$ and Koornwinder-Macdonald polynomials specialized at $t\sp {k+1}q\sp {r-1}=1$,
Compos. Math. 141 (2005), no. 6, 1589--1601.

\bibitem[Ka2]{Ka2}
M. Kasatani,
Subrepresentations in the polynomial representation of the double affine Hecke algebra of type ${\rm GL}\sb n$ at $t\sp {k+1}q\sp {r-1}=1$,
Int. Math. Res. Not. 2005, no. 28, 1717--1742.

\bibitem[KaSh]{KaSh}
M. Kasatani, K. Shigechi,
(in preparation.)

\bibitem[KaTa]{KaTa}
M. Kasatani, Y. Takeyama,
The quantum Knizhnik-Zamolodchikov equation and non-symmetric Macdonald polynomials,
Funkcialaj Ekvacioj, no. 50 (2007), 491--509.

\bibitem[Ma]{Ma}
I. Macdonald,
Affine Hecke Algebras and Orthogonal Polynomials,
Cambridge University Press, 2003.

\bibitem[NUKW]{NUKW}
A. Nishino, H. Ujino, Y. Komori, M. Wadati,
Rodrigues formulas for the nonsymmetric multivariable polynomials
 associated with the $BC\sb N$-type root system,
Nuclear Phys. B 571 (2000), no. 3, 632--648.

\bibitem[No]{No}
M. Noumi,
Macdonald-Koornwinder polynomials and affine Hecke rings,
S\=urikaisekikenky\=usho K\=oky\=uroku No. 919 (1995), 44--55 (in Japanese).

\bibitem[Ob]{Ob}
A. Oblomkov,
Double affine Hecke algebras of rank 1 and affine cubic surfaces,
Int. Math. Res. Not. 2004, no. 18, 877--912.

\bibitem[ObSt]{ObSt}
A. Oblomkov, E.Stoica,
Finite dimensional representations of double affine Hecke algebra of rank 1,
arXiv:math/0409256.

\bibitem[Sa]{Sa}
S. Sahi,
Nonsymmetric Koornwinder polynomials and duality,
Ann. of Math. (2) 150 (1999), no. 1, 267--282.

\bibitem[St]{St}
J.V. Stokman,
Koornwinder polynomials and affine Hecke algebras,
Internat. Math. Res. Notices 2000, no. 19, 1005--1042.


\end{thebibliography}
\end{document}